\documentclass[12pt]{article}

\input{epsf.sty}
\usepackage[dvips]{graphicx}
\usepackage{latexsym,amsmath,amsfonts,amscd,amsthm}
\usepackage{epsfig}
\usepackage{changebar}

\def\softd{{\leavevmode\setbox1=\hbox{d}%
          \hbox to 1.05\wd1{d\kern-0.4ex{\char039}\hss}}}%cstocs
\begin{document}
\baselineskip=2pc

%=============  title  =========================
\begin{center}

{\bf On the advantage of well-balanced schemes for moving-water
equilibria of the shallow water equations
}

\end{center}

\vspace{.02in}

\centerline{Yulong Xing\footnote{Computer Science and Mathematics
Division, Oak Ridge National Laboratory, Oak Ridge, TN 37831 and
Department of Mathematics, University of Tennessee, Knoxville, TN
37996. E-mail: xingy@math.utk.edu}, Chi-Wang Shu\footnote{Division
of Applied Mathematics, Brown University, Providence, RI 02912.
E-mail: shu@dam.brown.edu. Research supported by
AFOSR grant FA9550-09-1-0126 and NSF grant DMS-0809086.} and
Sebastian Noelle\footnote{Institute for Geometry and Applied
Mathematics, RWTH Aachen University, D-52056 Aachen, Germany. E-mail:
noelle@igpm.rwth-aachen.de. Research supported by DFG grant
GK 775.}}

\vspace{.2in}

%\centerline{January 2010}

\vspace{.2in}

\begin{abstract}

This note aims at demonstrating the advantage of moving-water
well-balanced schemes over still-water well-balanced schemes for
the shallow water equations.  We concentrate on numerical examples
with solutions near a moving-water equilibrium.  For such
examples, still-water well-balanced methods are not capable of
capturing the small perturbations of the moving-water equilibrium
and may generate significant spurious oscillations, unless an
extremely refined mesh is used. On the other hand, moving-water
well-balanced methods perform well in these tests.   The numerical
examples in this note clearly demonstrate the importance of
utilizing moving-water well-balanced methods for solutions near a
moving-water equilibrium.

\end{abstract}


\newpage

\section{Introduction}
\label{sec1} \setcounter{equation}{0} \setcounter{figure}{0}
\setcounter{table}{0}

The main objective of this note is to demonstrate the advantage of
schemes which are well balanced for moving water equilibrium,
over those which are only well balanced for still water, for small
perturbations of the moving water equilibrium for
the shallow water equations.

Well-balanced schemes refer to those schemes which have zero
truncation error for certain steady-state solutions of hyperbolic
balance laws
$$
U_t + f(U)_x = g(U,x) .
$$
That is, for certain non-trivial functions $V(x)$
satisfying
$$
f(V)_x = g(V,x),
$$
the well balanced schemes do not have any truncation error.
Note that in general $V(x)$ is unknown and is not a polynomial,
thus to require a zero truncation error is a difficult task.

For shallow water equations, there are typically two classes of
well-balanced schemes, namely those which are well balanced
for still-water equilibrium, where the velocity is zero, and
those which are well balanced for general moving-water
equilibrium.  A scheme which belongs to the latter class
is more difficult to construct.  We refer to
\cite{NXS2007} for a recently developed high-order accurate
finite volume weighted essentially non-oscillatory (WENO)
scheme which is well balanced for general moving-water
equilibrium.  Discussions on the history of the development
of well balanced schemes for shallow water equations,
with an extensive reference list, can be found in \cite{NXS2010}.

Since it is much more difficult to construct well-balanced
schemes for moving-water equilibrium than for still water,
a natural question is whether the former has any advantage
over the latter.  However, to our best knowledge, examples to show this advantage
seem to be absent in the literature.  In the
Castro Urdiales meeting for which this special issue is
dedicated to, there was an intense discussion
on this issue but no such examples were provided in the
discussion.

The purpose of this note is therefore to provide carefully
selected numerical examples to demonstrate the
advantage of well-balanced
schemes for moving-water equilibrium over
well-balanced schemes for still water.
These examples have solutions which contain small perturbations
from a moving-water equilibrium.  A well-balanced scheme
for still water would generate spurious oscillations at the
level of the truncation error of the scheme, until the
mesh is extremely refined.  On the other hand, a
well-balanced scheme for the moving-water equilibrium could
resolve these small perturbations on much coarser meshes.
These examples clearly demonstrate the advantage of
well-balanced schemes for moving-water equilibrium.
We shall use the high order finite volume WENO
schemes in \cite{XS2006b}, which are well balanced
for still water, and those in \cite{NXS2007}, which
are well balanced for moving water, in the numerical tests.

\section{Shallow water equations and well-balanced methods}
\label{sec:lakeatrest} \setcounter{equation}{0}
\setcounter{figure}{0} \setcounter{table}{0}

In one space dimension, the shallow water equations take the form
\begin{equation}  \label{main1}
\left\{%
\begin{array}{l}
\displaystyle h_t+(hu)_x=0 \\
\displaystyle (hu)_t+\left( hu^2+\frac{1}{2}gh^2 \right)_x=-ghb_x ,%
\end{array}
\right.
\end{equation}
where $h$ denotes the water height, $u$ is the velocity of the
fluid, $b$ represents the bottom topography and $g$ is the
gravitational constant. The still-water (also referred as lake at
rest) steady state is given by
\begin{equation}  \label{stationary}
u=0 \qquad \text{and} \qquad h+b=const,
\end{equation}
which represents a still flat water surface. The general
moving-water steady state is given by
\begin{equation}  \label{steady}
m=hu=const \qquad \text{and} \qquad E=\frac{1}{2}u^2+g(h+b)=const.
\end{equation}
We introduce the notation $U=(h,\,hu)^{T}$ for the conservative
variables, with the superscript $T$ denoting the transpose, and
$V=(m,\,E)^T$ for the equilibrium variables. We refer to
\cite{NXS2007} for the variable transformation between
$U=U(V,x)$ and $V=V(U,x)$.

The computational domain is discretized into cells
$I_{i}=[x_{i-\frac{1}{2}},x_{i+\frac{1}{2}}]$, $i=1, \cdots ,N$.
We denote the size of the $i$-th cell by $\triangle x_{i}$ and the
center of the cell by $x_i=\left(
x_{i-\frac{1}{2}}+x_{i+\frac{1}{2}} \right)/2$. Let
$\bar{U}(x_i,t)=\int_{I_i}U(x,t) \, dx/\triangle x_i$ denote the
cell average of $U( \cdot ,t)$ over the cell $I_i$. The
computational variables are $\bar{U}^n_i$, which approximate the
cell averages $\bar{U}(x_i,t)$ at the time $t=t^n$. For the ease
of presentation, we denote the shallow water equations
(\ref{main1}) by
\begin{equation*}
U_{t}+f(U)_{x}=g(h,b).
\end{equation*}

\subsection{Methods preserving still-water equilibrium}

Many high-order well-balanced methods for the still-water steady
state solution (\ref{stationary}) have been
developed in the literature for the shallow water equations,
see for example \cite{NXS2010} for a list of references.
In this note, we will use the fifth-order finite volume WENO
methods based on the separation of the source terms developed in
Xing and Shu \cite{XS2006b}.
In this subsection, we briefly review this well-balanced scheme in one
dimension.  We refer to \cite{XS2006b} for further details.

We first reconstruct the point values $U^{\pm}_{i+\frac{1}{2}}$ at the
cell interface from the given cell
averages $\bar{U}_i$ by the standard WENO reconstruction
procedure, which can be eventually written out as
\begin{equation}
\label{add1}
 U_{i+\frac{1}{2}}^+ = \sum_{k=-r+1}^{r} w_k \bar{U}_{i+k} \equiv
S^+_{\bar{U}} (\bar{U})_i , \qquad
 U_{i+\frac{1}{2}}^- = \sum_{k=-r}^{r-1} \tilde{w}_k \bar{U}_{i+k} \equiv
S^-_{\bar{U}} (\bar{U})_i
\end{equation}
where $r=3$ for the fifth order WENO approximation and the
coefficients $w_k$ and $\tilde{w}_k$ depend nonlinearly on
$\bar{U}$.
% changes by Sebastian
Now we freeze the cell averages $\bar{U}$ and the corresponding coefficients
$w_k$ and $\tilde{w}_k$. With these fixed coefficients, we apply the {\em linear}
operators $W \mapsto S^{\pm}_{\bar{U}}(W)$ to $W=(b,0)^{T}$
to compute the reconstructed values $b^{\pm}_{i+\frac{1}{2}}$.
Assume $U$ is the still-water steady state solution satisfying
(\ref{stationary}), we clearly have
\begin{equation}
h^{\pm}_{i+\frac{1}{2}} + b^{\pm}_{i+\frac{1}{2}} = constant.
\end{equation}

The main idea in constructing well-balanced methods is to
decompose the integral of the source term into a sum of several
terms, then compute each of them in a way consistent with the
approximation for the corresponding flux terms. We first rewrite
the shallow water equations (\ref{main1}) as
\begin{equation}
  \left\{\begin{array} {l}
         \displaystyle   h_t+(hu)_x=0 \\
         \displaystyle
         (hu)_t+ \left( hu^2+\frac{1}{2}gh^2 \right)_x=
\left(\frac12gb^2\right)_x-g(h+b)b_x,
         \end{array}
   \right.
\end{equation}
for which the semi-discrete numerical method takes the form
\begin{equation}
  \label{main3}
 \triangle x_i \frac{d\bar{U}(x_i,t)}{dt} =
        -(\hat{F}_{i+\frac{1}{2}}-\hat{F}_{i-\frac{1}{2}})
        +G(U,x),
\end{equation}
where $\hat{F}$ is a numerical flux, such as the
Lax-Friedrichs flux
\begin{equation}
 \label{lf2}
 \hat{F}_{i+\frac{1}{2}} =
   \frac{1}{2}\left[ f(U^-_{i+\frac12})+f(U^+_{i+\frac12})
                   - \alpha\left(\left( \begin{array}{c}
                                         h+b \\ hu
                                        \end{array} \right)^+_{i+\frac{1}{2}}
                               - \left( \begin{array}{c}
                                         h+b \\ hu
                                        \end{array} \right)^-_{i+\frac{1}{2}}
                           \right) \right],
\end{equation}
with $\alpha=\max_U|u\pm\sqrt{gh}|$. The first component of the
source term approximation $G(U,x)$ is zero, and the second
component is given by
\begin{equation}
\frac{1}{2}g\left(\widehat{b_{i+\frac12}^2}-\widehat{b_{i-\frac12}^2}\right)
- g\overline{(h+b)}_i (\hat{b}_{i+\frac12}- \hat{b}_{i-\frac12})
-\int_{I_i} g\left(h+b-\overline{(h+b)}_i\right) b_x dx,
\end{equation}
with
$\widehat{b^2_{i+\frac12}}=\frac12\left((b^+_{i+\frac12})^2
+(b^-_{i+\frac12})^2\right)$
and
$\hat{b}_{i+\frac12}=\frac12(b^+_{i+\frac12}+b^-_{i+\frac12})$.

\subsection{Methods preserving moving-water equilibrium}

We have developed high-order well-balanced finite volume WENO
schemes for the moving-water equilibrium (\ref{steady}) for the
shallow water equations in \cite{NXS2007}. In this subsection, we
briefly review these methods and refer to \cite{NXS2007} for
further details.

At each time step, we first apply the usual WENO reconstruction
procedure to the cell averages $\bar{U}_i$, and obtain
$U^{\pm}_{i+\frac{1}{2}}$, hence $V^{\pm}_{i+\frac{1}{2}}$.

Given the cell averages $\bar{U_i}$ and a bottom function $b(x)$,
we choose local reference values $\bar{V_i}$ of the equilibrium
variables. These are defined implicitly by the requirement that
\begin{equation}
    \label{eq:Vb_i}
    \frac{1}{\Delta x_i}\int_{I_i}U(\bar{V}_i,x)dx = \bar{U}_i.
\end{equation}
Relation \eqref{eq:Vb_i} chooses $\bar{V}_i$ as the unique (see
the paragraphs preceding \cite[Def.3.2]{NXS2007}) local
equilibrium such that the corresponding conserved variables
$U(\bar{V}_i,b(x))$ have the same cell average $\bar{U}_i$ as the
numerical data. It is proven in \cite[Def.3.2]{NXS2007} that, if
the data $U(x)$ and $b(x)$ are in local equilibrium
($V(U(x),x)\equiv \bar{V}$ for all cells $I_i$), the reference
equilibrium states $\bar{V}_i$ computed via \eqref{eq:Vb_i}
coincide with the true local steady state $\bar{V}$. The
reconstruction is completed by limiting the reconstruction
$V^{\pm}_{i+\frac{1}{2}}$ with respect to the reference values
$\bar{V}_i$ (see \cite[(3.18)]{NXS2007}), to obtain the limited
values $\tilde{V}^{\pm}_{i+\frac{1}{2}}$ and $\tilde{V}_{i}$.

The fourth-order well-balanced scheme is given by
\begin{equation}
  \label{main4}
 \triangle x_i \frac{d\bar{U}(x_i,t)}{dt} =
        -(\hat{F}(\hat U_{i+\frac12}^-,\hat U_{i+\frac12}^+)
         -\hat{F}(\hat U_{i-\frac12}^-,\hat U_{i-\frac12}^+))
        +s_i,
\end{equation}
where the function $\hat{F}(\cdot,\cdot)$ is a conservative, Lipschitz
continuous numerical flux, and
\begin{equation}\label{eq:U_B} \hat
U_{i+\frac12}^\pm = U(\tilde V_{i+\frac12}^\pm, \hat
b_{i+\frac12}), \qquad \hat
b_{i+\frac12}=\min(b_{i+\frac12}^{-},b_{i+\frac12}^{+}).
\end{equation}

The total source term $s_i$ is given by
\begin{equation}\label{eq:s_wb1}
s_i := \frac{\displaystyle 4S_2-S_1}{3} +f(\hat U_{i-\frac12}^+)
-f(\tilde U_{i-\frac12}^+) +f(\tilde U_{i+\frac12}^-) -f(\hat
U_{i+\frac12}^-),
\end{equation}
where $\tilde U_{i-\frac12}^{\pm} = U(\tilde V_{i+\frac12}^\pm,
b_{i+\frac12}^\pm)$. The extrapolated interior source term
$(4S_2-S_1)/3$ is defined by
\begin{eqnarray}\label{eq:s_wb2} S_1
&=& s_i^{int}(\tilde U^+_{i-\frac{1}{2}}, \tilde
U^-_{i+\frac{1}{2}},
 b^+_{i-\frac{1}{2}}, b^-_{i+\frac{1}{2}}) \\
\label{eq:s_wb3} S_2 &=& \left( s_i^{int}(\tilde U^+_{i-\frac12},
\tilde U_i, b^+_{i-\frac{1}{2}}, b_i) +s_i^{int}(\tilde U_i,
\tilde U^-_{i+\frac12}, b_i, b^-_{i+\frac12}) \right)
\end{eqnarray}
and the well-balanced quadrature of the source term $s_i^{int}$ is
given by
\begin{equation}\label{eq:s_wb4}
s_i^{int}(U_L,U_R,b_L,b_R) = -\frac12g(h_L+h_R)(b_R-b_L)+\hat
s_i^{int}
\end{equation}
where $\hat s_i^{int}$ is given by \cite[(3.63)-(3.65)]{NXS2007}.

\section{Numerical examples}
\label{sec:numerical} \setcounter{equation}{0}
\setcounter{figure}{0} \setcounter{table}{0}

In this section we present numerical results of both moving-water
well-balanced methods and still-water well-balanced methods
presented in Section \ref{sec:lakeatrest}, for a carefully selected set
of test examples
in one dimension. Comparison of these results is provided as a
demonstration of the advantage of moving-water well-balanced
methods. In all the examples, time discretization is by the
classical third-order total variation diminishing (TVD)
Runge-Kutta method \cite{eno1}, and the
CFL number is taken as 0.6. The gravitation constant $g$ is taken
as 9.812$m/s^2$.

\subsection{Perturbation of a moving water equilibrium}\label{sec:n1}

The following test cases are chosen to demonstrate the capability
of these schemes for computations on the perturbation of steady
state solutions.

The bottom topography is given by:
\begin{equation}
\label{7.3.1} b(x) = \left\{
      \begin{array} {l l}
       0.2-0.05(x-10)^2 & \mbox{if } 8\leq x\leq 12, \\
       0 & \mbox{otherwise,}
      \end{array}
      \right.
\end{equation}
in the computational domain $[0,25]$. Three steady states,
subcritical or transcritical flow with or without a steady shock
will be investigated.

Our initial conditions are given by imposing a small perturbation
of size $0.05$ on the height of these steady states in the
interval $[5.75,\,6.25]$. Theoretically, this disturbance should
split into two waves, propagating to the left and right
respectively. Note that in \cite{NXS2007}, we have shown the
numerical results of the moving-water well-balanced methods with a
perturbation size $0.01$, which demonstrate that these small
perturbations are well captured. Here we provide the numerical
results of both the moving-water well-balanced methods and
still-water well-balanced methods, and demonstrate the different
behavior of these two methods.

a): {\em Subcritical flow:}

The initial condition is given by:
\begin{equation}\label{eq:subcritical}
E = 22.06605, \qquad m = 4.42,
\end{equation}
together with the boundary condition that the discharge $m$=4.42
is imposed at upstream and the water height $h$=2 is imposed at
downstream when the flow is subcritical.

\begin{figure}
\centerline{
\includegraphics[width=3.2in]{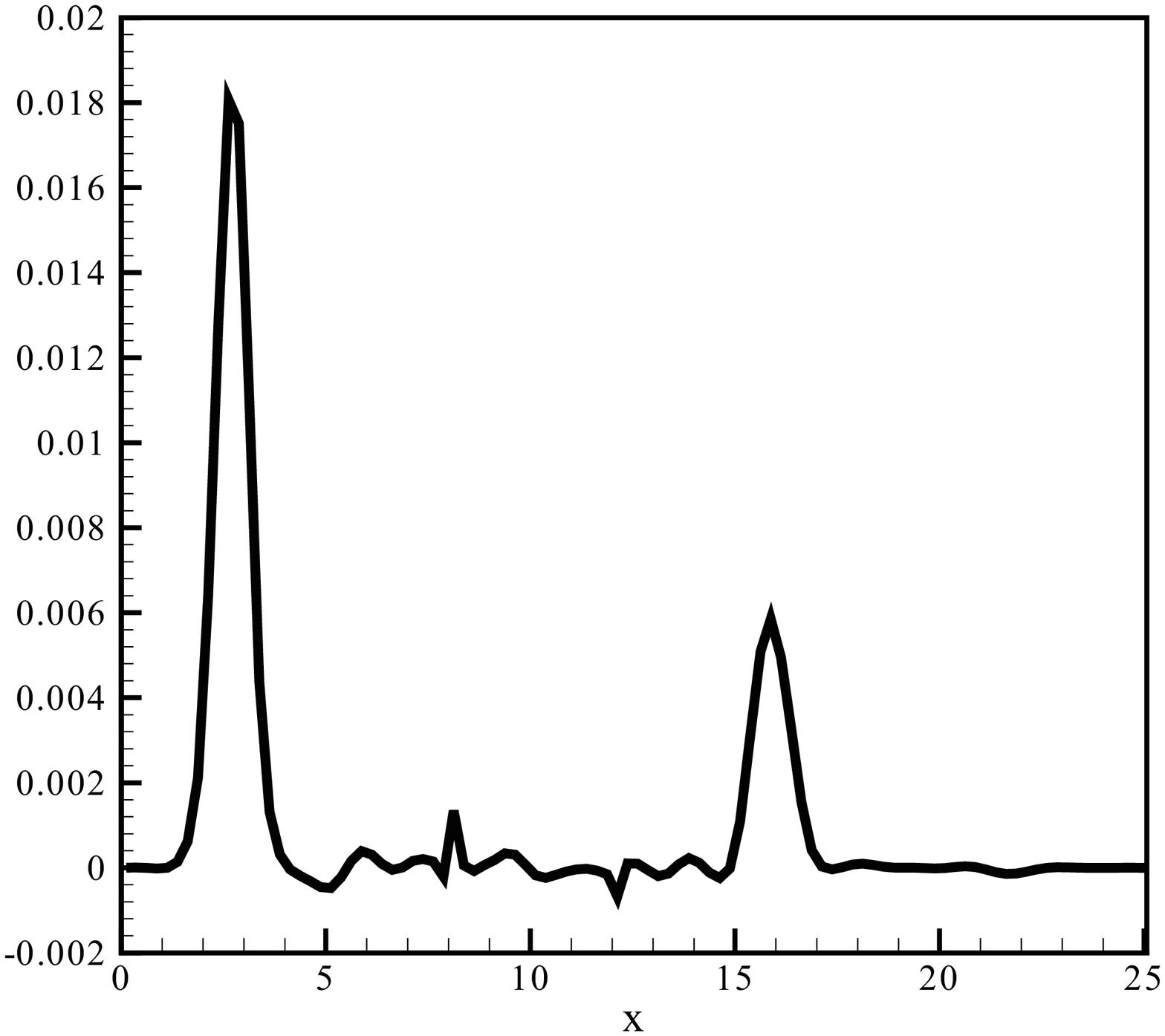}
\hspace{0.05in}
\includegraphics[width=3.2in]{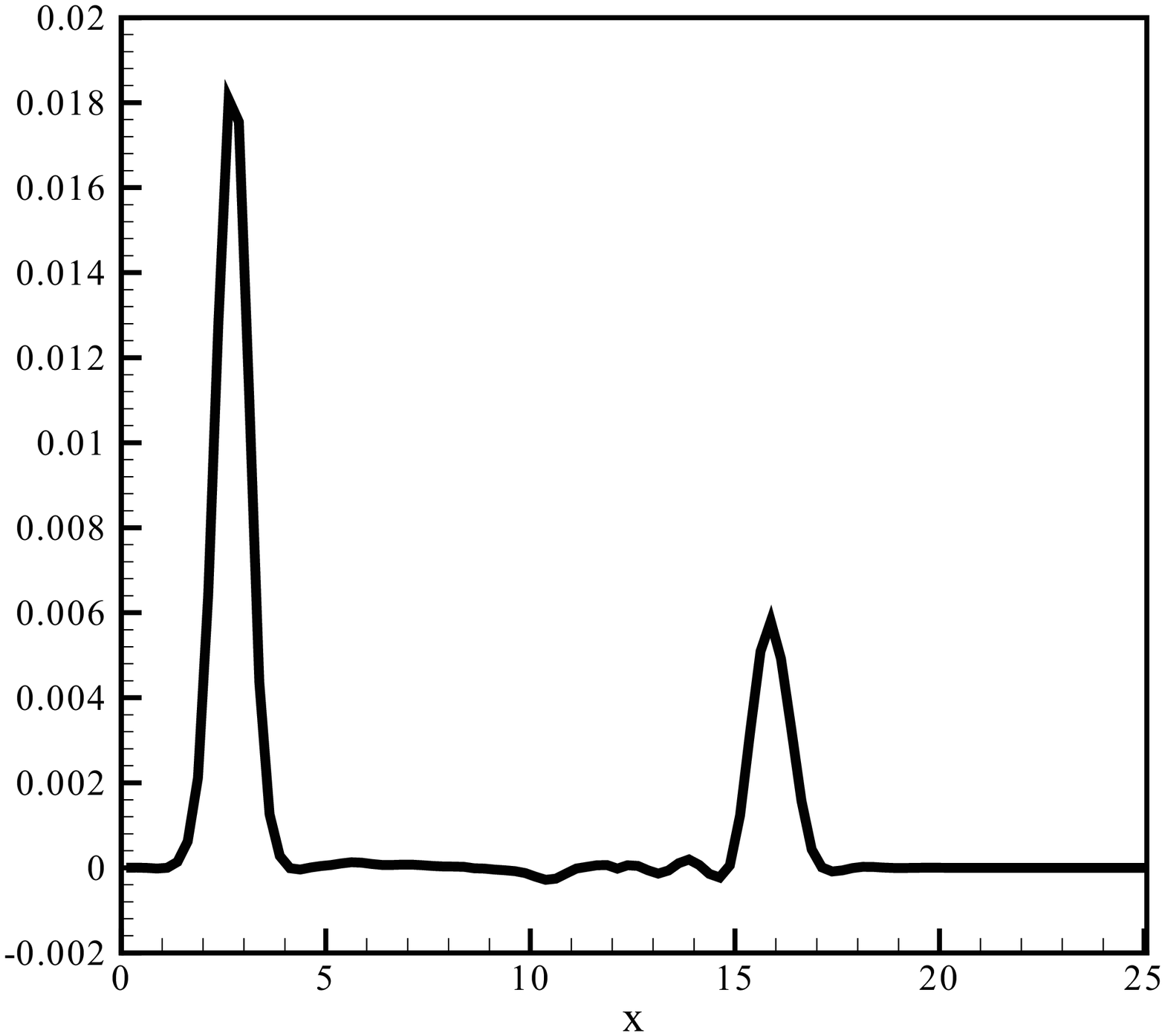}
} \vspace{-0.15in} \caption{The difference between the height $h$
at time $t=1.5$ and the background moving steady state water
height (\ref{eq:subcritical}), when 100 uniform cells are
employed. An initial perturbation of size $0.05$ is imposed
between $[5.75,\,6.25]$. Left: result based on still-water
well-balanced scheme. Right: result based on moving-water
well-balanced scheme.} \label{f1}
\end{figure}

\begin{figure}
\centerline{
\includegraphics[width=3.2in]{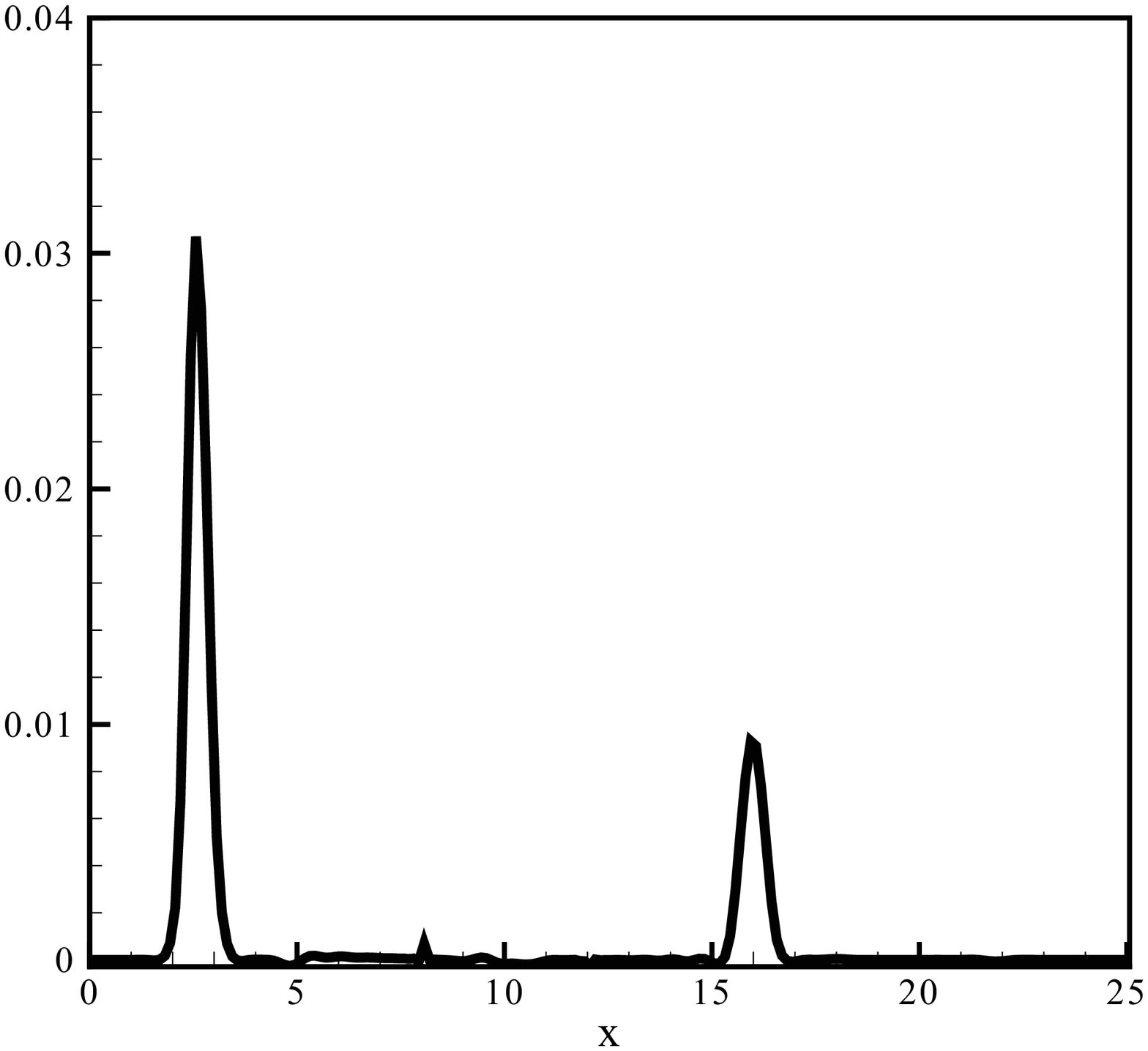}
\hspace{0.05in}
\includegraphics[width=3.2in]{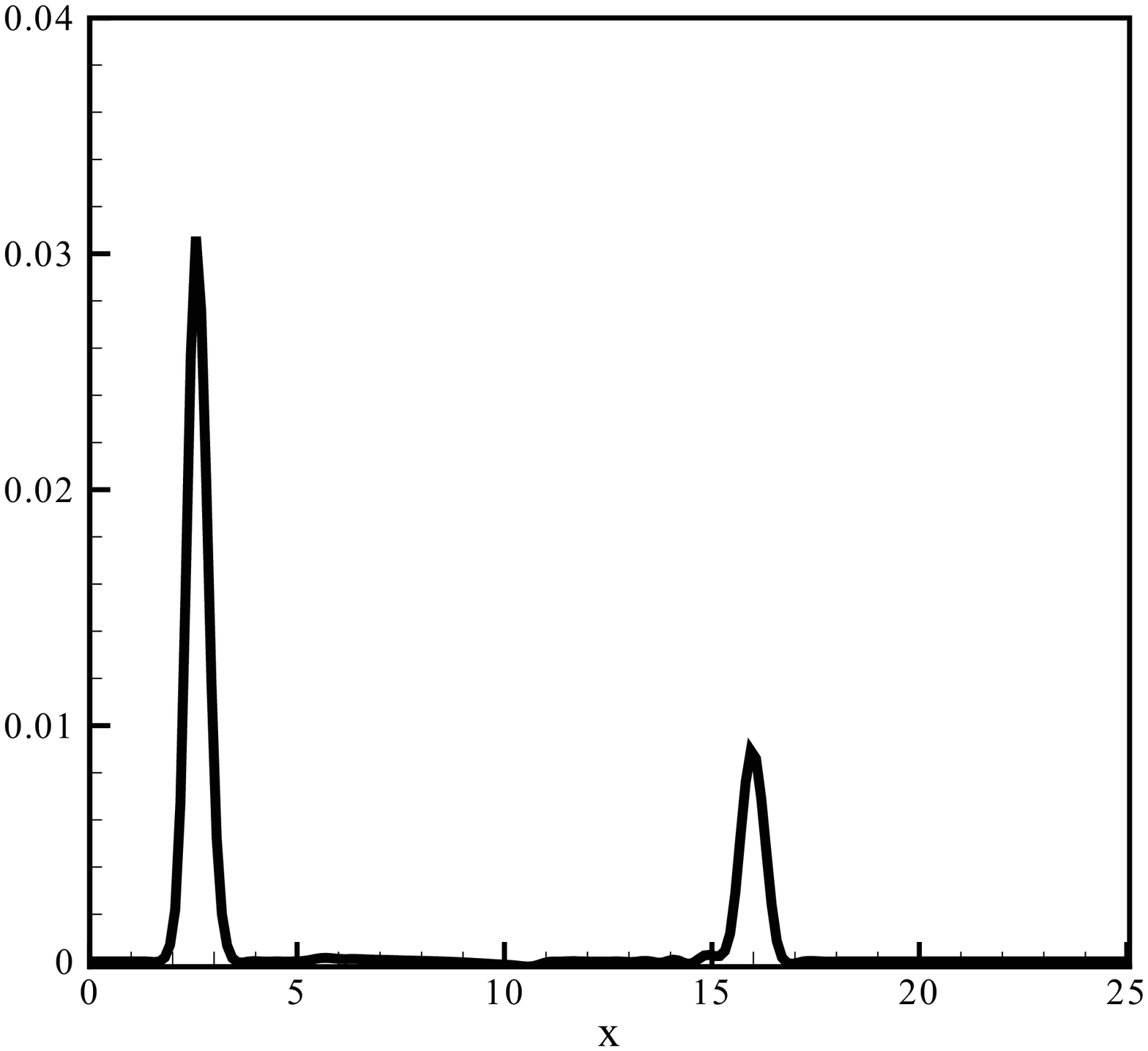}
} \vspace{-0.15in} \caption{Same as in Figure \ref{f1}, but 1000
uniform cells are employed.} \label{f2}
\end{figure}

We run the tests with 100 uniform cells until $t=1.5$, when
the downstream-traveling water pulse has already passed the bump. The
differences between the water height $h$ at that time and the
background moving water state, are shown in Figure \ref{f1}.
The result of the still-water well-balanced method is plotted on the
left and that of the moving-water well-balanced method
is shown on the right.  The difference can be easily observed, as
there are significantly more spurious oscillations generated by the still-water
well-balanced method. As we refine the mesh to $1000$ cells, the
results are shown in Figure \ref{f2}. Because $\Delta x$ is now
very small, the perturbation is relatively big in comparison
with truncation errors of the schemes and can be well
captured by both methods.

b): {\em Transcritical flow without a shock:}

The initial condition is given by:
\begin{equation} \label{eq:transhock}
E = \frac{1.53^2}{2\times0.66^2} + 9.812\times0.66, \qquad m =
1.53,
\end{equation}
together with the boundary condition that the discharge $m$=1.53
is imposed at upstream and the water height $h$=0.66 is imposed at
downstream when the flow is subcritical.

\begin{figure}
\centerline{
\includegraphics[width=3.2in]{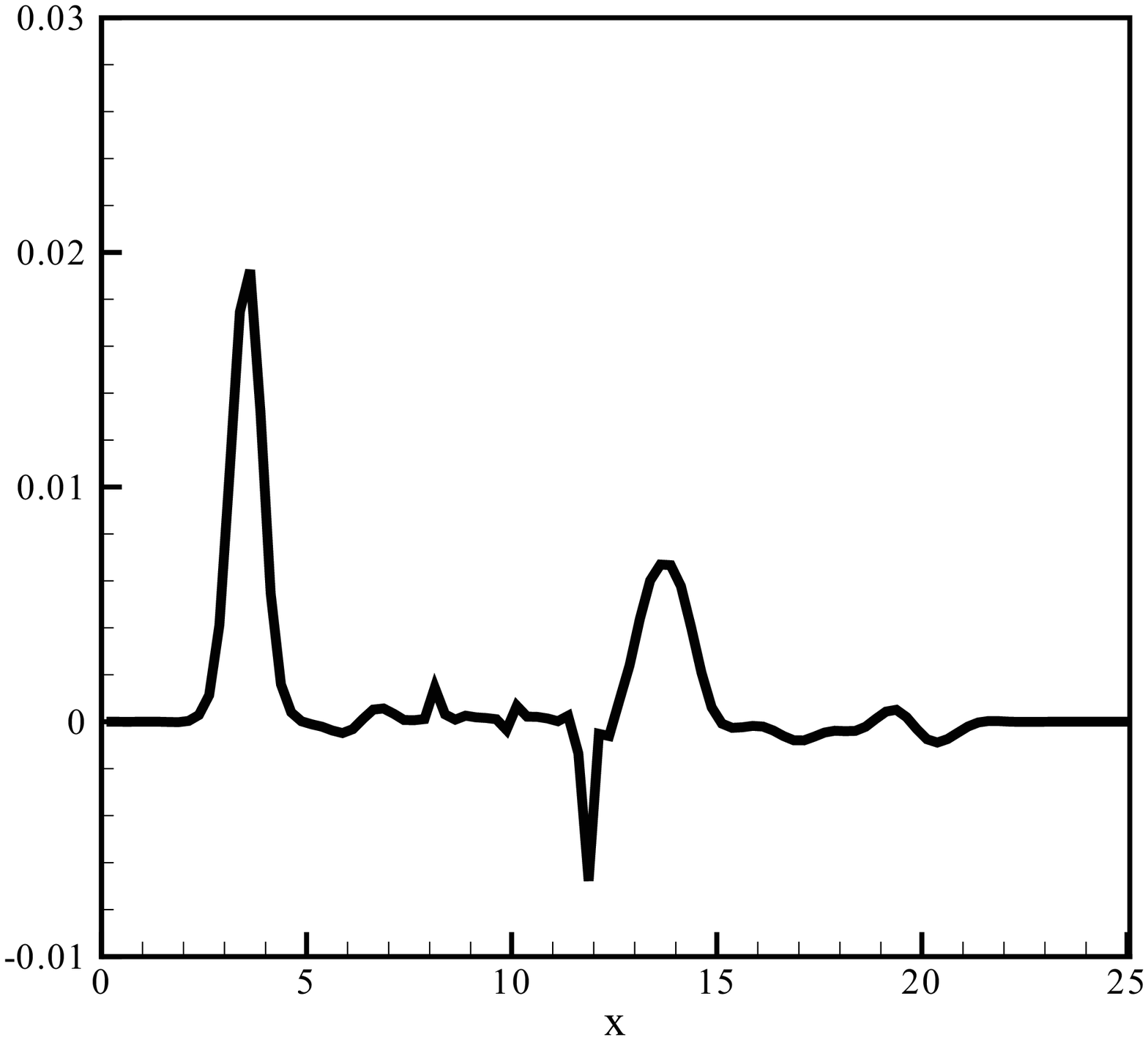}
\hspace{0.05in}
\includegraphics[width=3.2in]{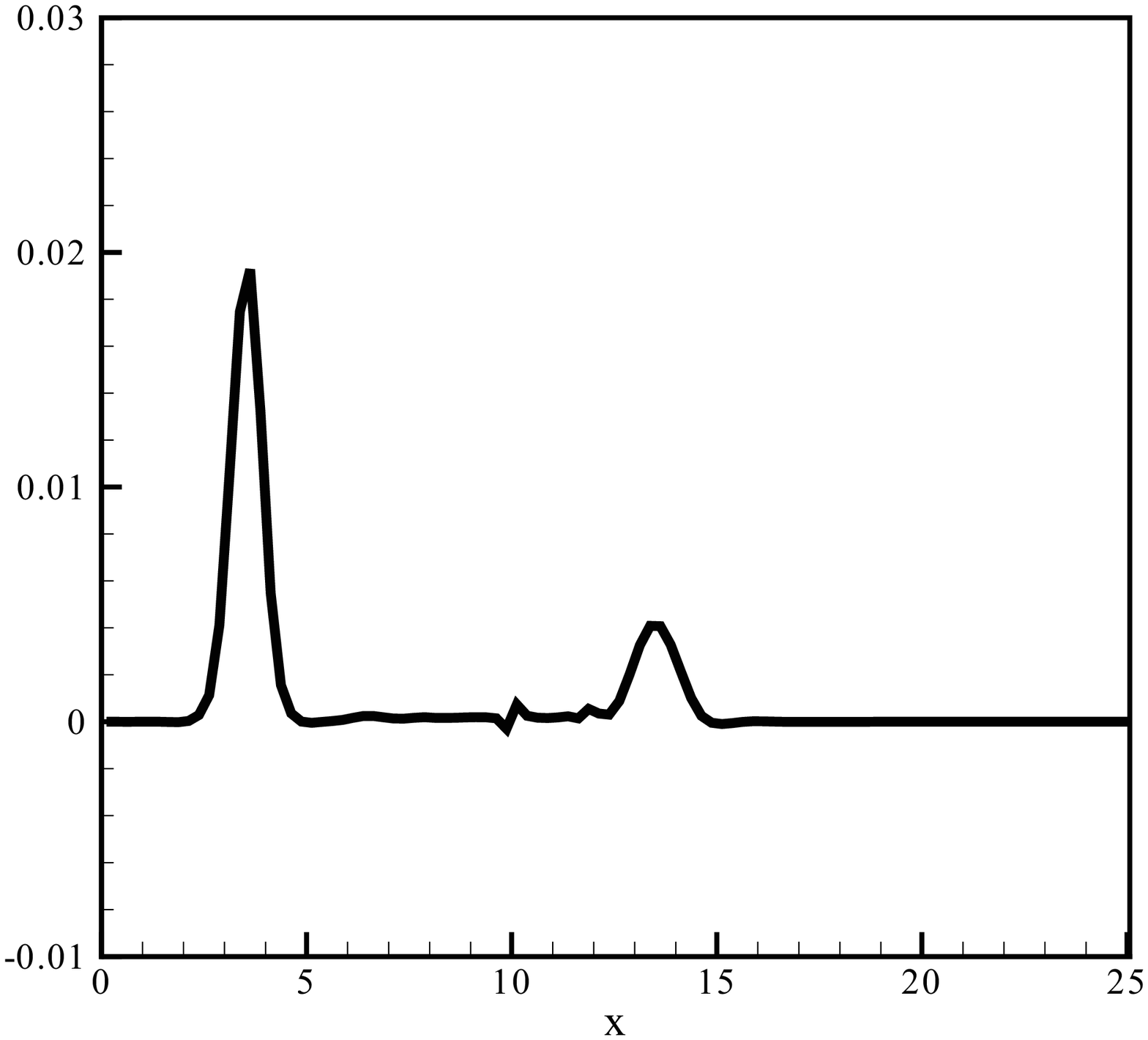}
} \vspace{-0.15in} \caption{The difference between the height $h$
at time $t=1.5$ and the background moving steady state water
height (\ref{eq:transhock}), when 100 uniform cells are employed.
An initial perturbation of size $0.05$ is imposed between
$[5.75,\,6.25]$. Left: result based on still-water well-balanced
scheme. Right: result based on moving-water well-balanced scheme.}
\label{f3}
\end{figure}

\begin{figure}
\centerline{
\includegraphics[width=3.2in]{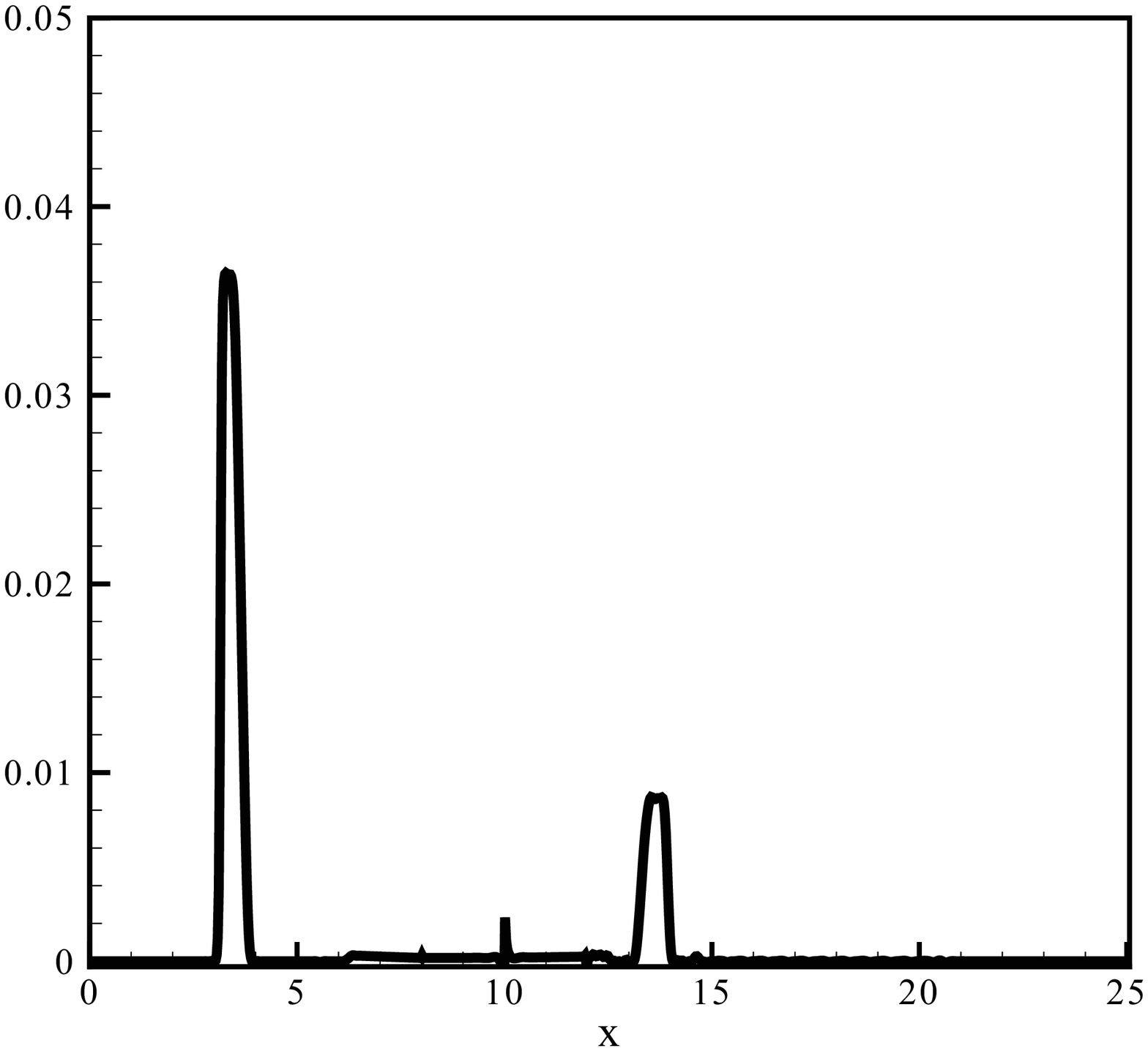}
\hspace{0.05in}
\includegraphics[width=3.2in]{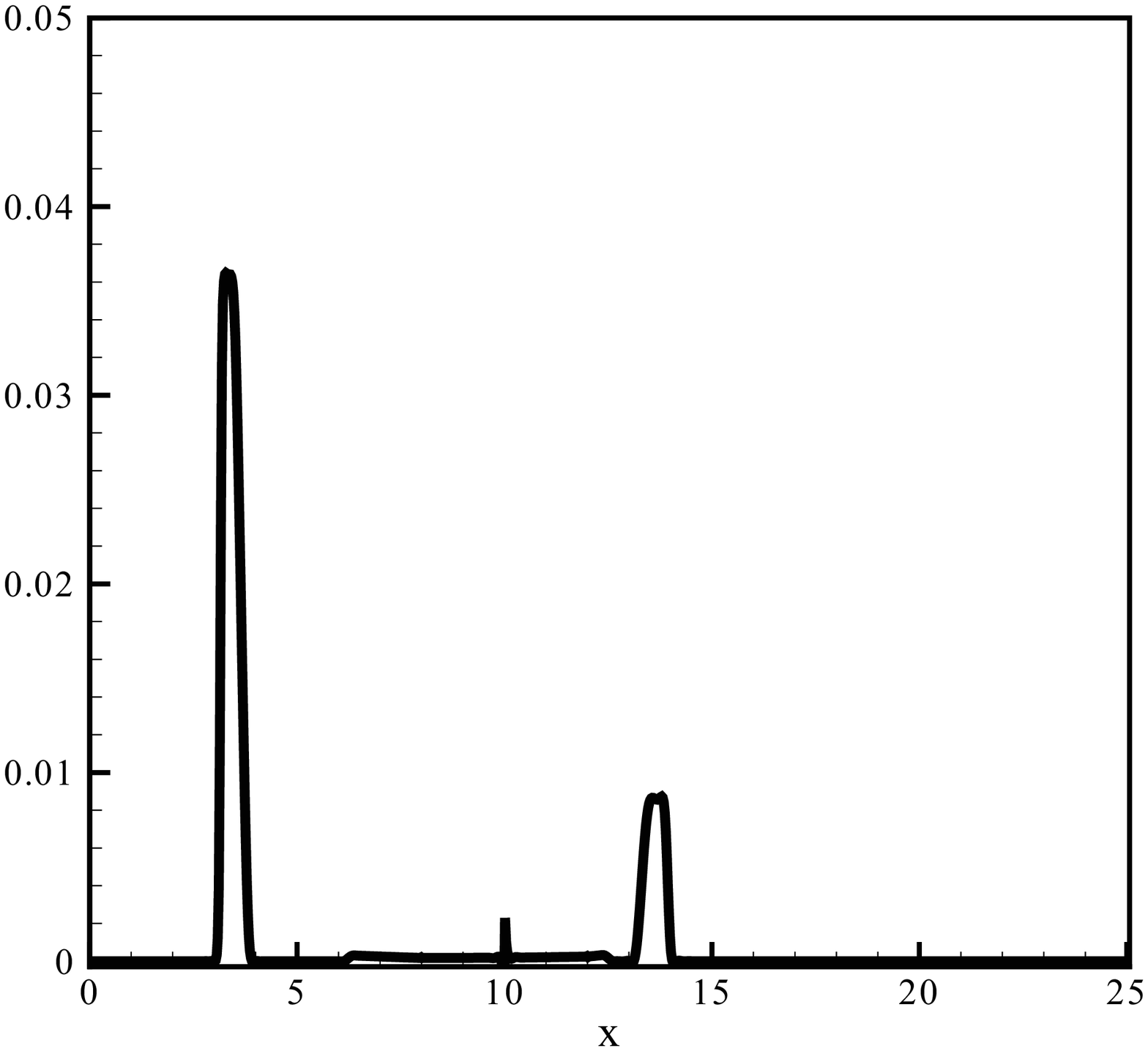}
} \vspace{-0.15in} \caption{Same as in Figure \ref{f3}, but 1000
uniform cells are employed.} \label{f4}
\end{figure}

We run the tests with 100 uniform cells until $t=1.5$, when
the downstream-traveling water pulse has already passed the bump. The
differences between the water height $h$ at that time and the
background moving water state, are shown in Figure \ref{f3}.
Again, the difference can be easily observed, as there are significantly more
spurious oscillations generated by the still-water well-balanced method.
As we refine the mesh to $1000$ cells, the results in
Figure \ref{f4} show resolved solutions for both methods.

c): {\em Transcritical flow with a shock:}

The initial condition is given by:
\begin{equation} \label{eq:transhock1}
E = \left\{
      \begin{array} {l l}
       \frac32(9.812 \times 0.18)^{\frac23}) + 9.812 \times 0.2

         & \mbox{if } x\leq 11.665504281554291   \\
       \frac{\displaystyle 0.18^2}{\displaystyle 2\times0.33^2}
+ 9.812\times0.33 & \mbox{otherwise}
      \end{array}
      \right.
 \qquad m = 0.18,
\end{equation}
together with the boundary condition that the discharge $m$=0.18
is imposed at upstream and the water height $h$=0.33 is imposed at
downstream.

\begin{figure}
\centerline{
\includegraphics[width=3.2in]{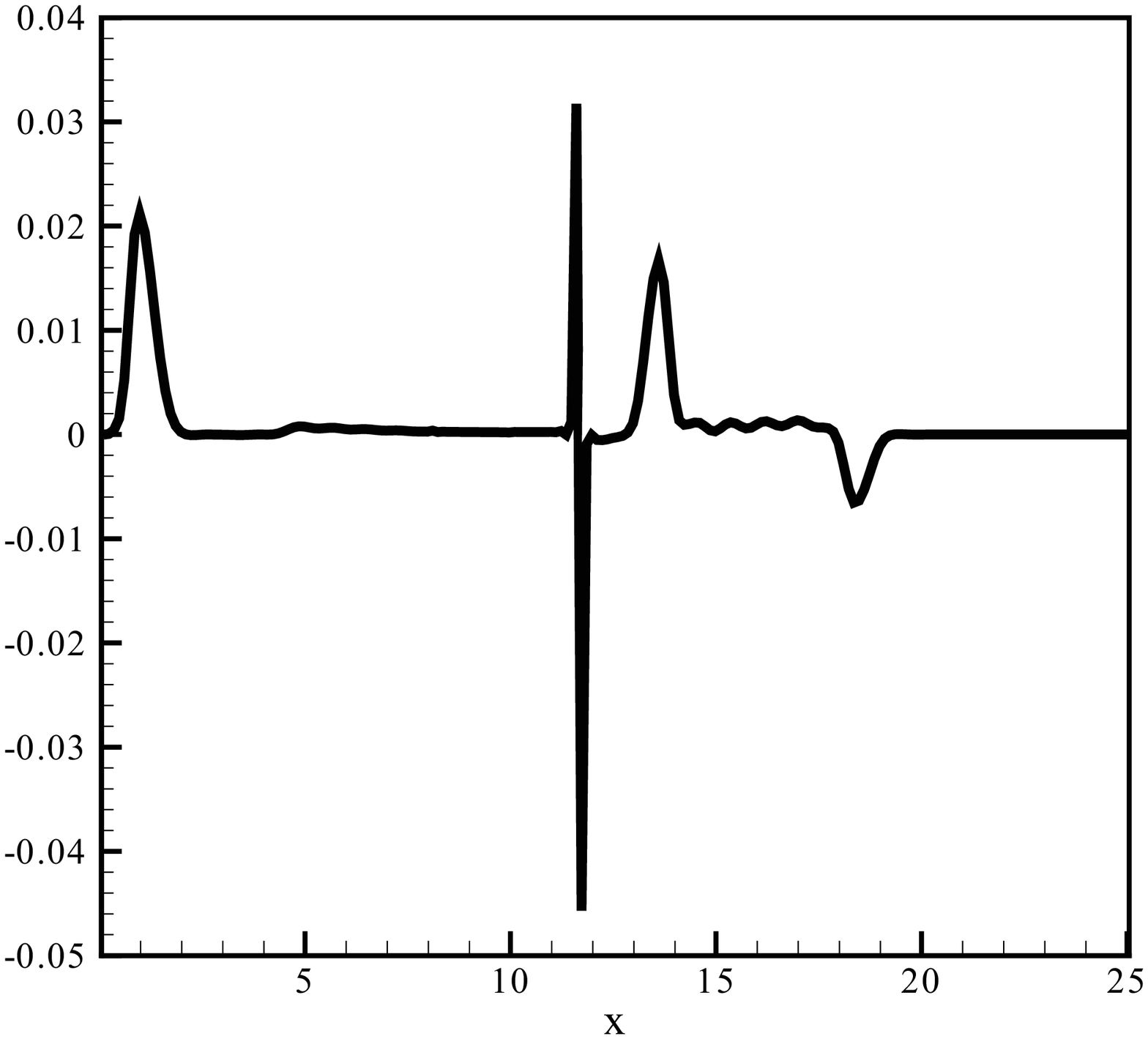}
\hspace{0.05in}
\includegraphics[width=3.2in]{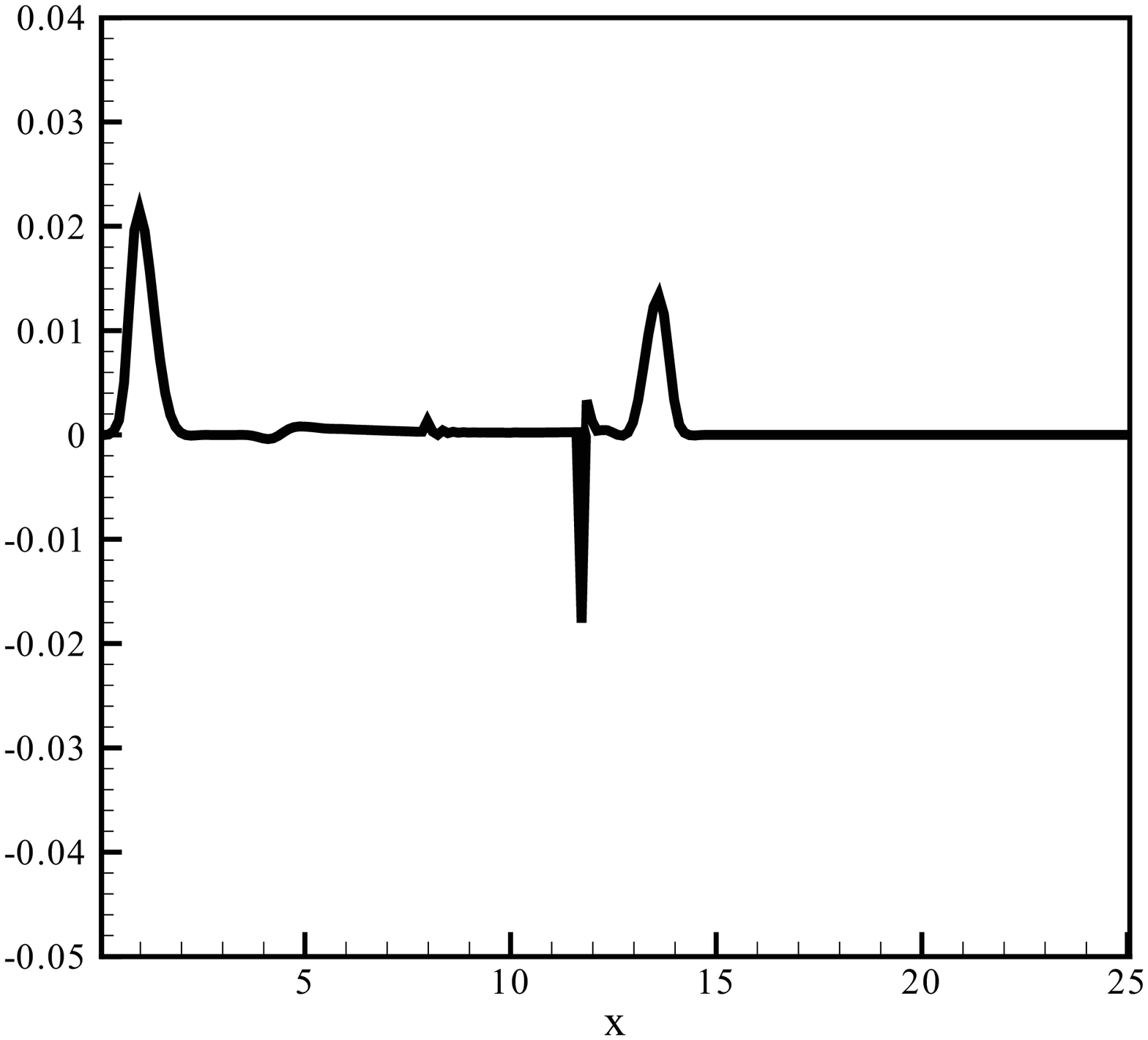}
} \vspace{-0.15in} \caption{The difference between the height $h$
at time $t=3$ and the background moving steady state water height
(\ref{eq:transhock1}), when 200 uniform cells are employed. An
initial perturbation of size $0.05$ is imposed between
$[5.75,\,6.25]$. Left: result based on still-water well-balanced
scheme. Right: result based on moving-water well-balanced scheme.}
\label{f5}
\end{figure}

\begin{figure}
\centerline{
\includegraphics[width=3.2in]{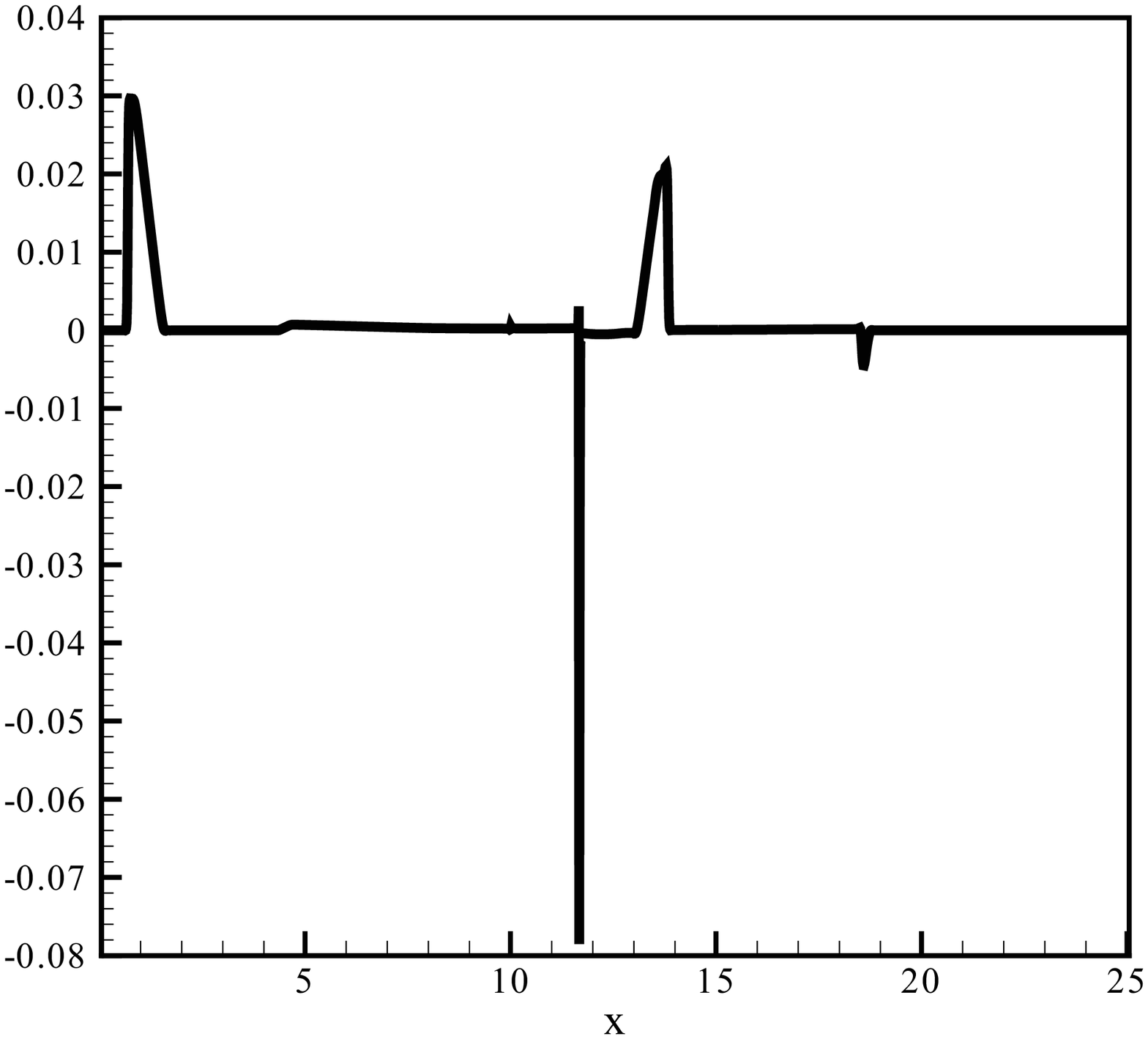}
\hspace{0.05in}
\includegraphics[width=3.2in]{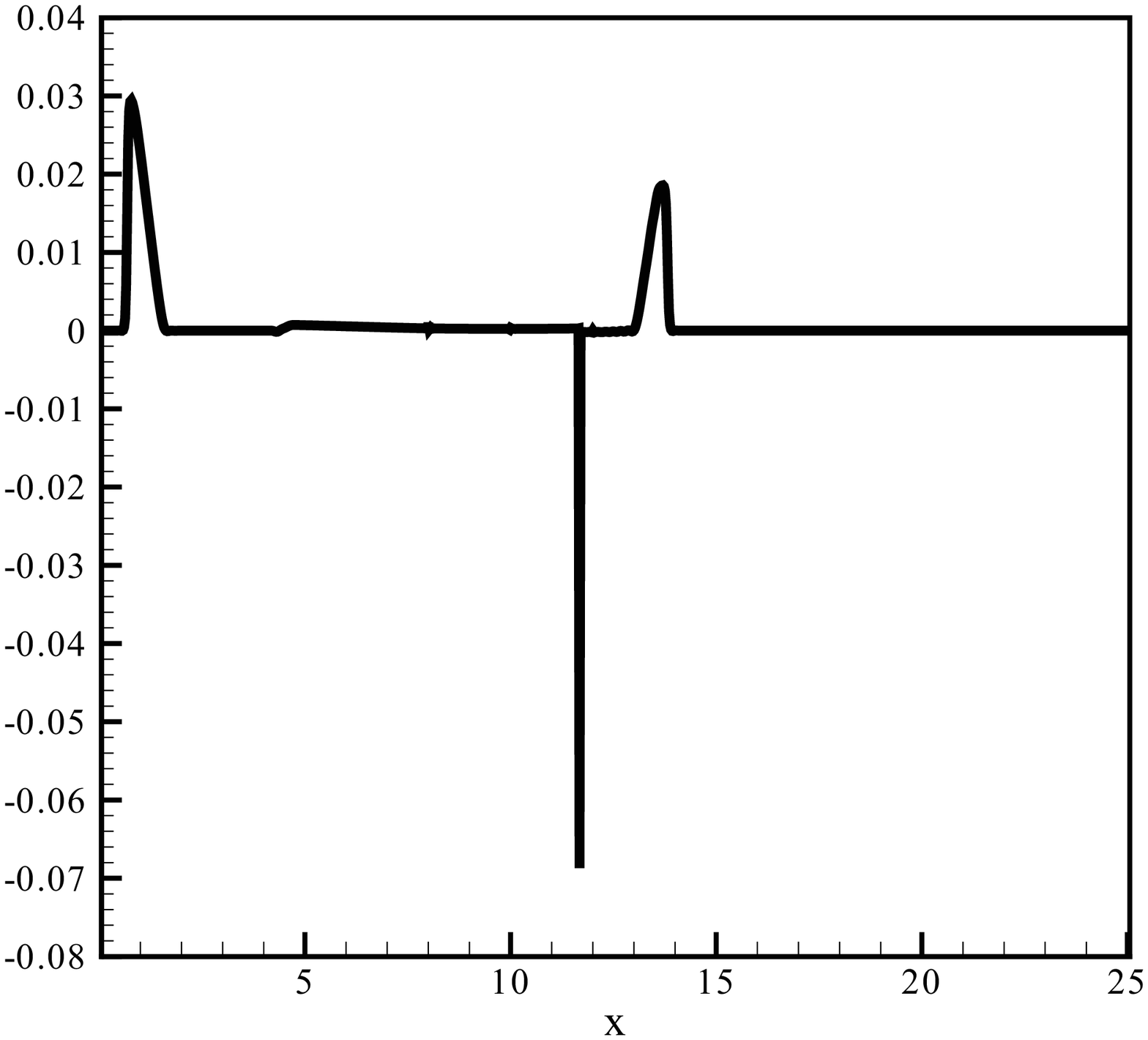}
} \vspace{-0.15in} \caption{Same as in Figure \ref{f5}, but 1000
uniform cells are employed.} \label{f6}
\end{figure}

We run the tests with 200 uniform cells until $t=3$. The
differences between the water height $h$ at that time and the
background moving water state, are shown in Figure \ref{f5}.
Once again, the difference can be easily observed, as there are significantly more
spurious structures generated by the still-water well-balanced method.
As we refine the mesh to $1000$ cells, the results are shown in
Figure \ref{f6}.  Even at this resolution there is still an advantage
using the moving-water well-balanced method.
Note that near
the point $x=11.7$ there is a big displacement for both methods,
as we recall that $11.7$ is the position where the shock is located.

\subsection{Perturbation with smaller magnitude}\label{sec:n2}

In this subsection, we utilize a smaller perturbation to these
tests and show that moving-water well-balanced methods demonstrate
more clearly their advantage in capturing such perturbations.

\begin{figure}
\centerline{
\includegraphics[width=3.2in]{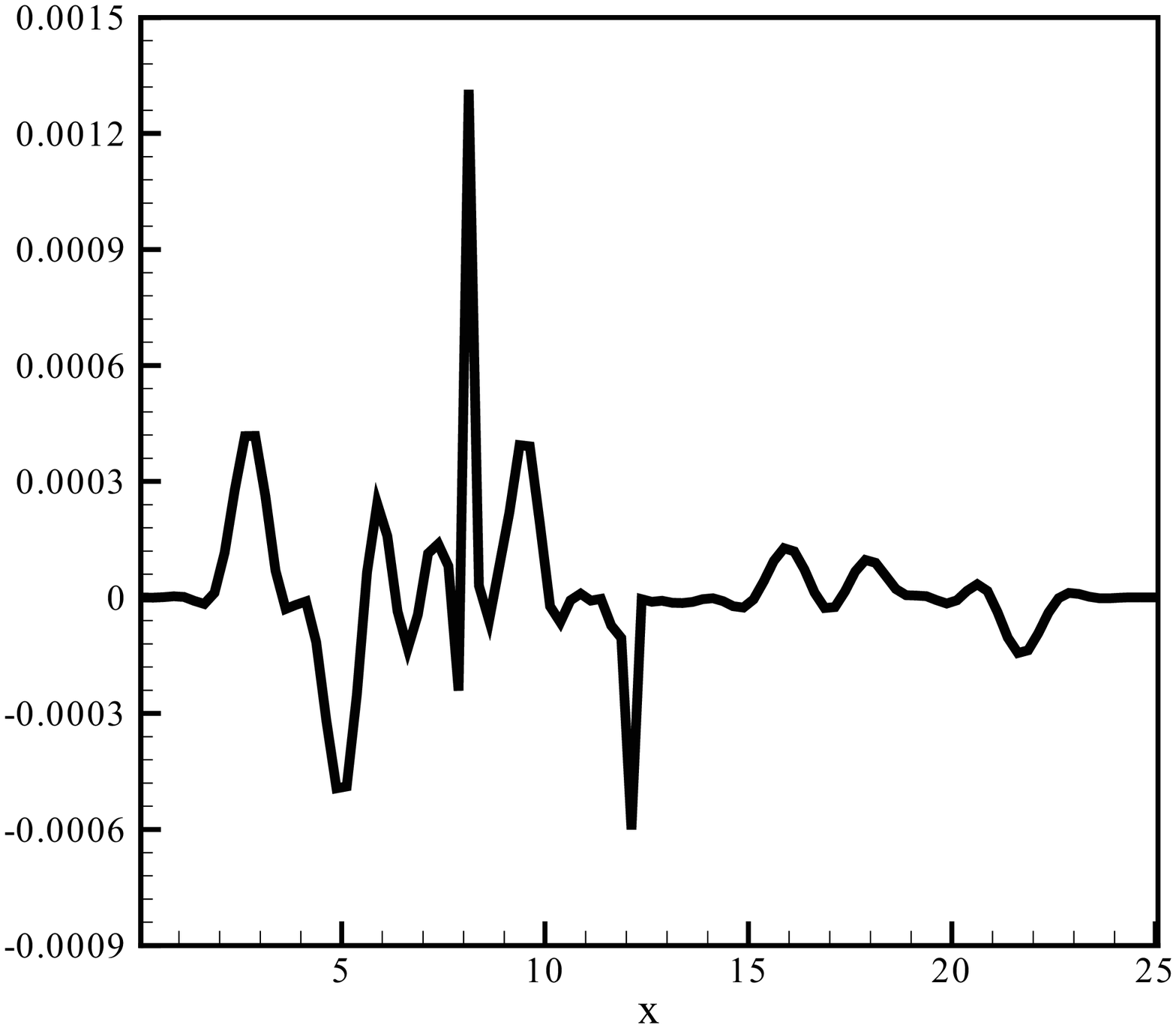}
\hspace{0.05in}
\includegraphics[width=3.2in]{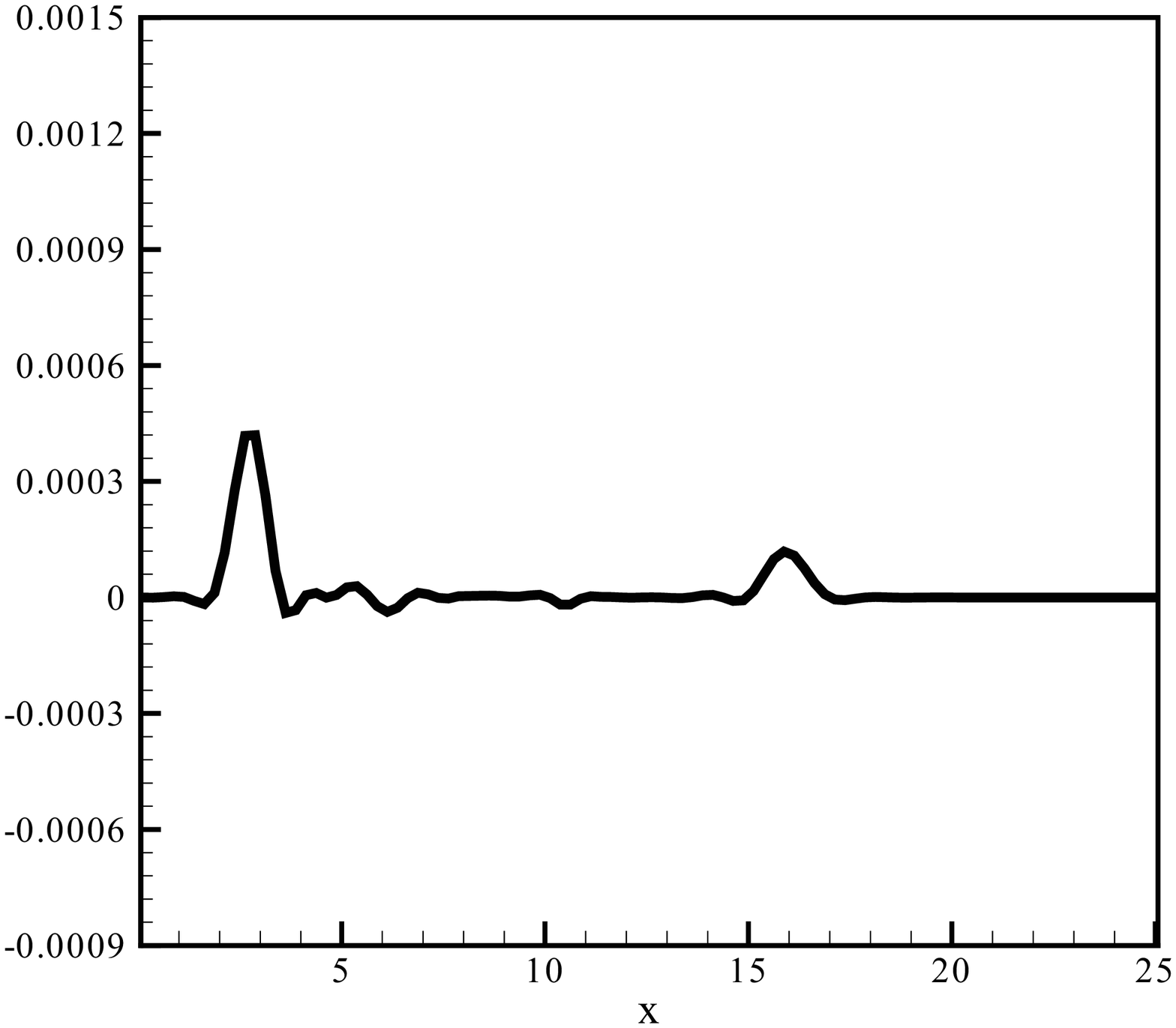}
} \vspace{-0.15in} \caption{The difference between the height $h$
at time $t=1.5$ and the background moving steady state water
height (\ref{eq:subcritical}), when 100 uniform cells are
employed. An initial perturbation of size $0.001$ is imposed
between $[5.75,\,6.25]$. Left: result based on still-water
well-balanced scheme. Right: result based on moving-water
well-balanced scheme.} \label{f11}
\end{figure}

\begin{figure}
\centerline{
\includegraphics[width=3.2in]{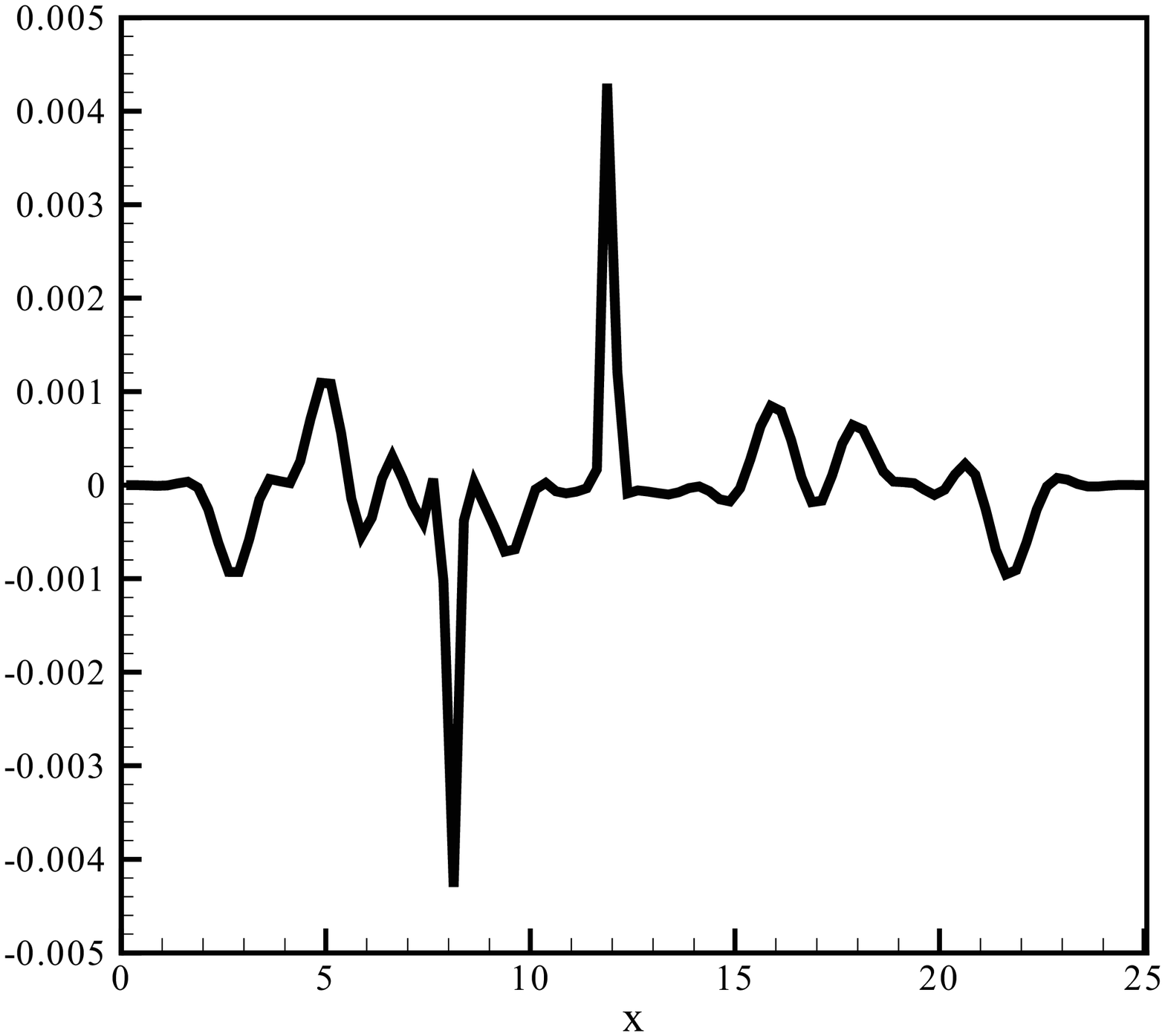}
\hspace{0.05in}
\includegraphics[width=3.2in]{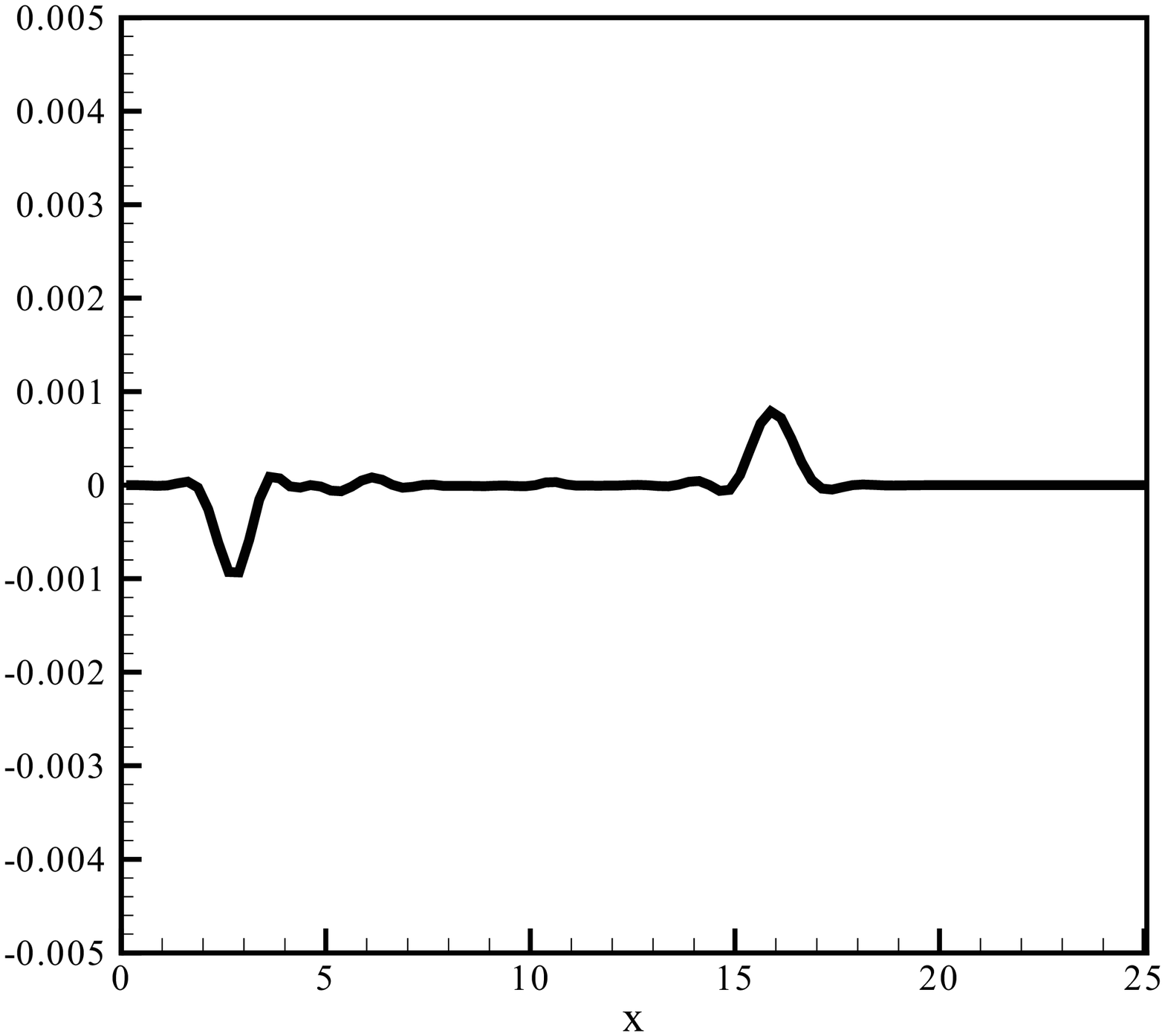}
} \vspace{-0.15in} \caption{The difference between the momentum
$hv$ at time $t=1.5$ and the background moving steady state water
height (\ref{eq:subcritical}), when 100 uniform cells are
employed. An initial perturbation of size $0.001$ is imposed
between $[5.75,\,6.25]$. Left: result based on still-water
well-balanced scheme. Right: result based on moving-water
well-balanced scheme.} \label{f21}
\end{figure}

\begin{figure}
\centerline{
\includegraphics[width=3.2in]{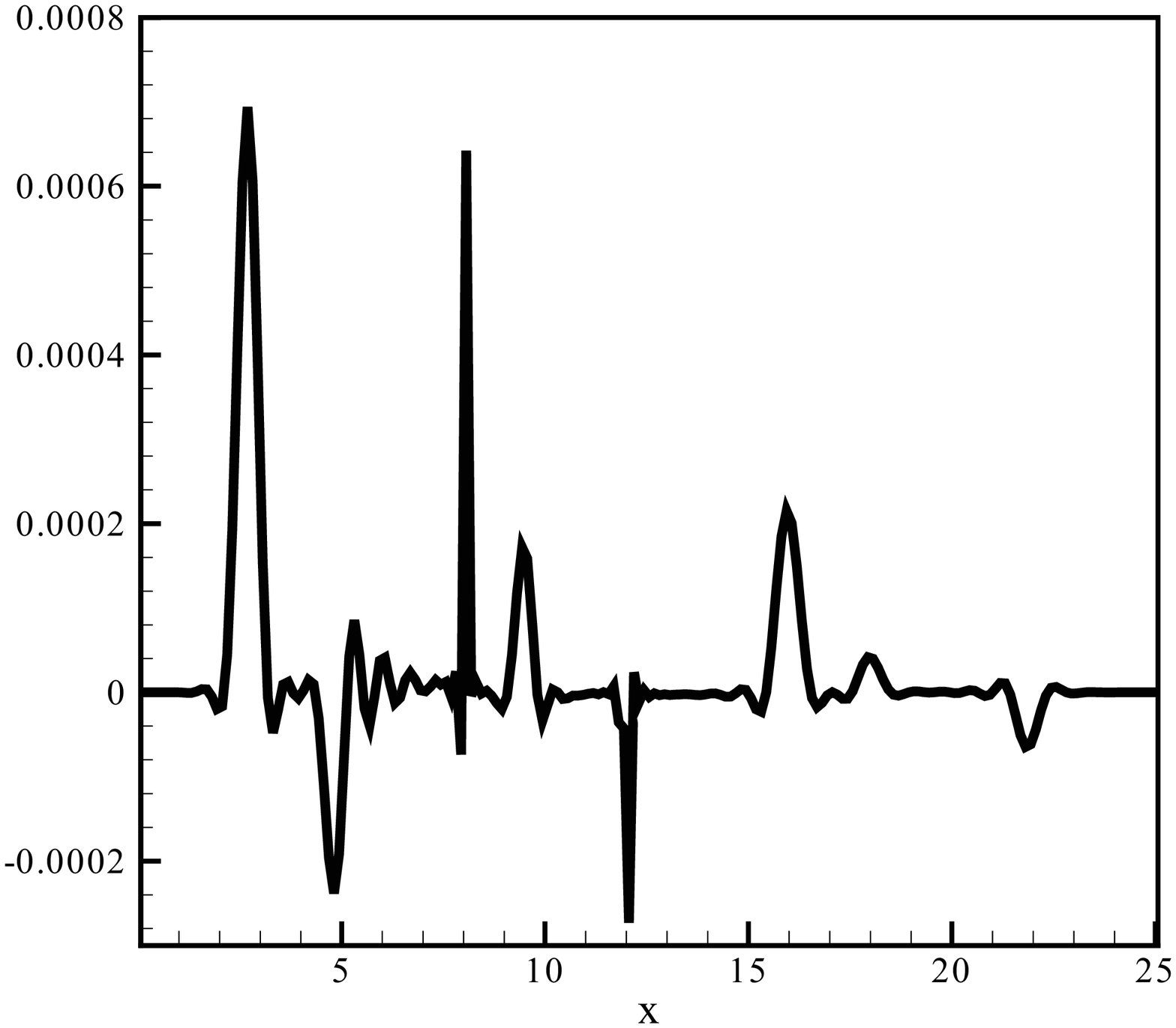}
\hspace{0.05in}
\includegraphics[width=3.2in]{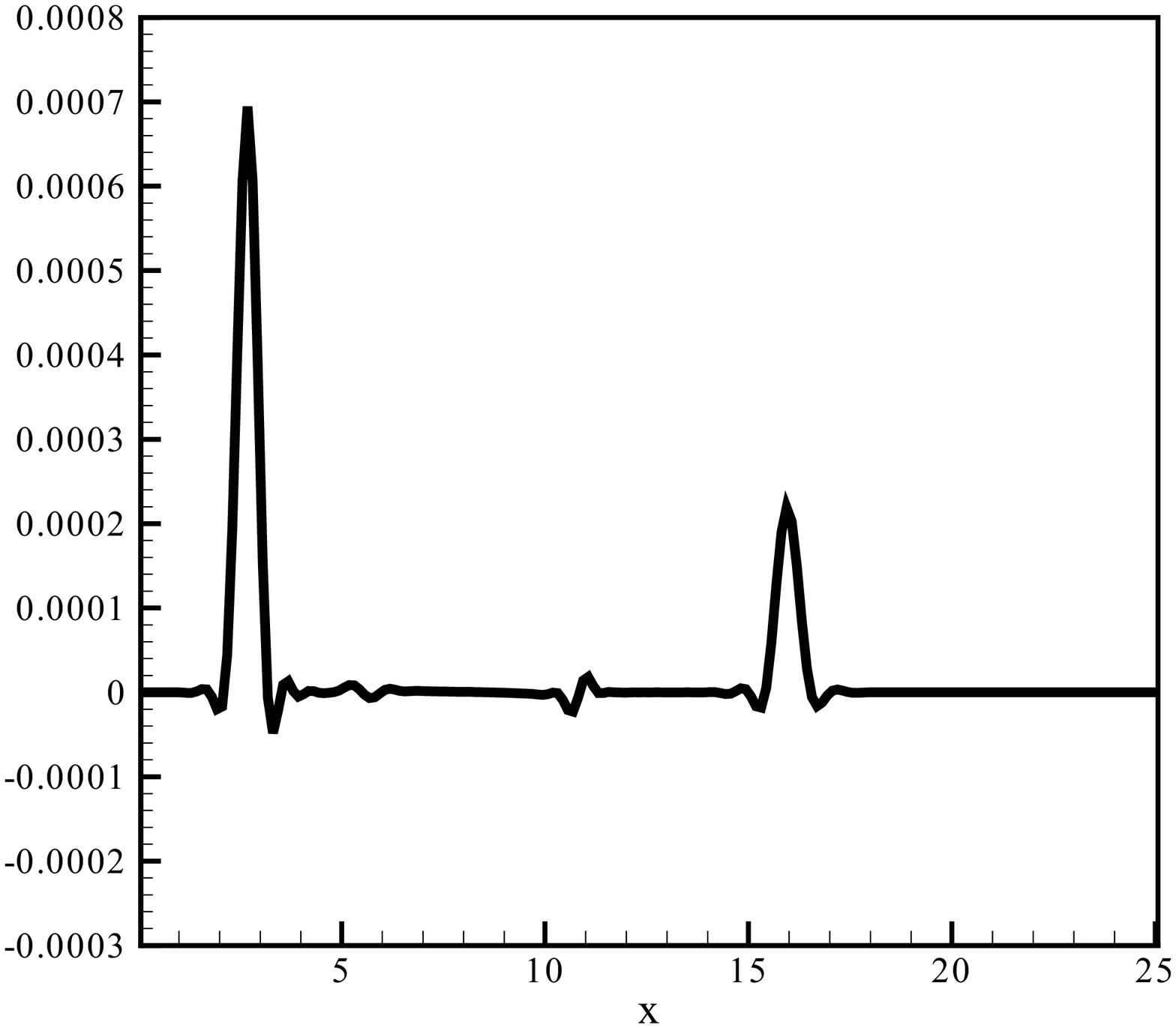}
} \vspace{-0.15in} \caption{Same as in Figure \ref{f11}, but 200
uniform cells are employed.} \label{f12}
\end{figure}

\begin{figure}
\centerline{
\includegraphics[width=3.2in]{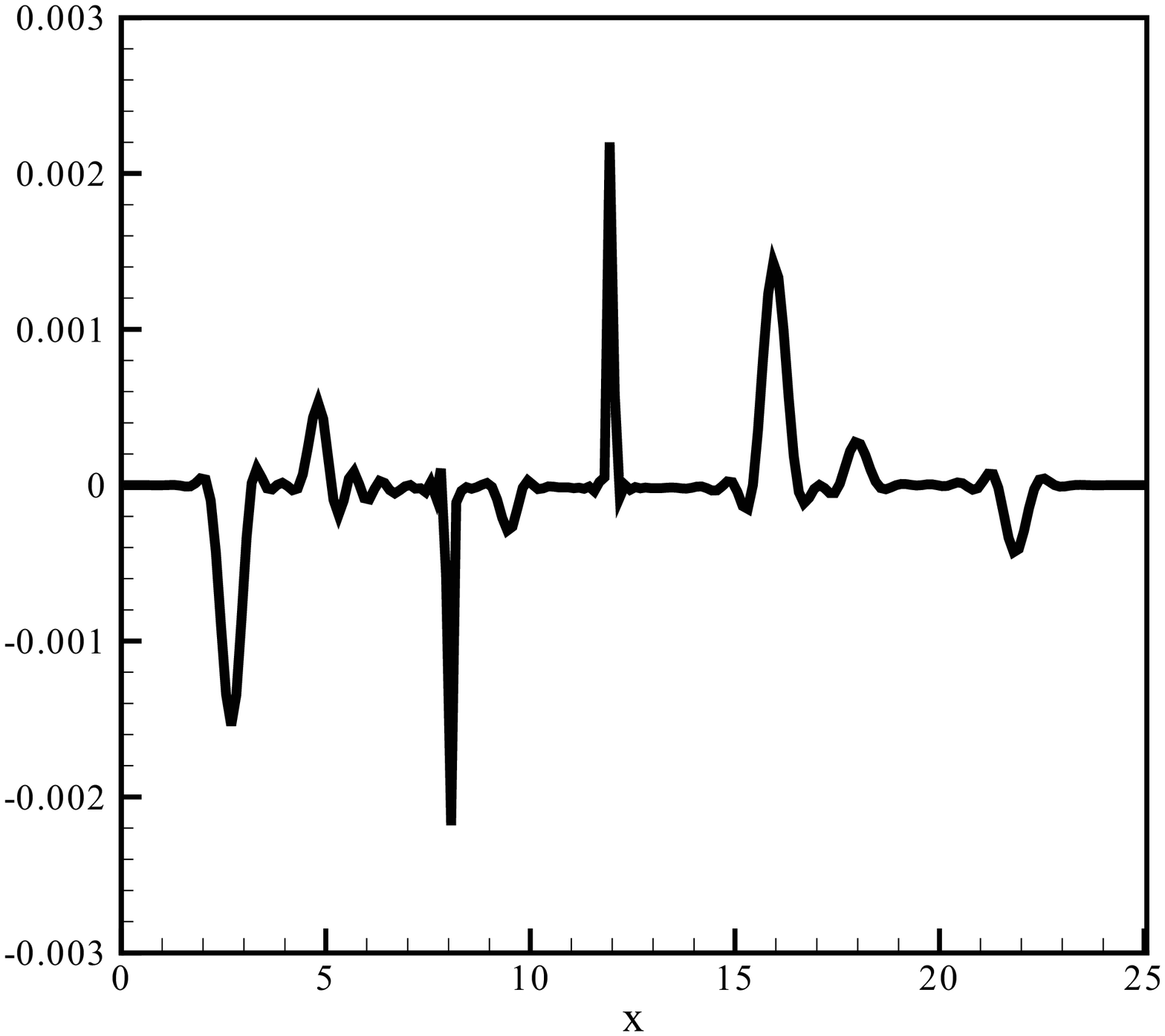}
\hspace{0.05in}
\includegraphics[width=3.2in]{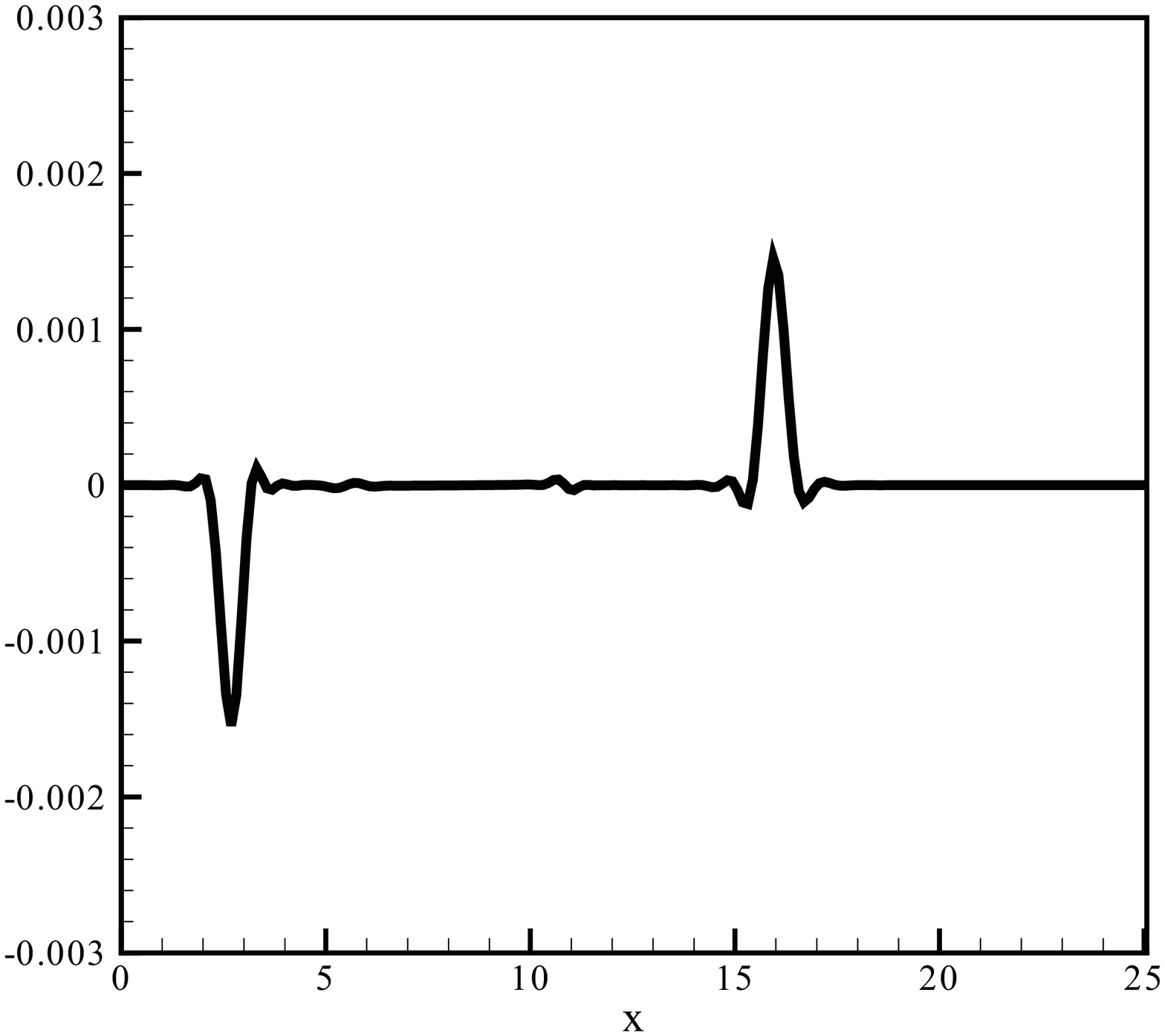}
} \vspace{-0.15in} \caption{Same as in Figure \ref{f21}, but 200
uniform cells are employed.} \label{f22}
\end{figure}

\begin{figure}
\centerline{
\includegraphics[width=3.2in]{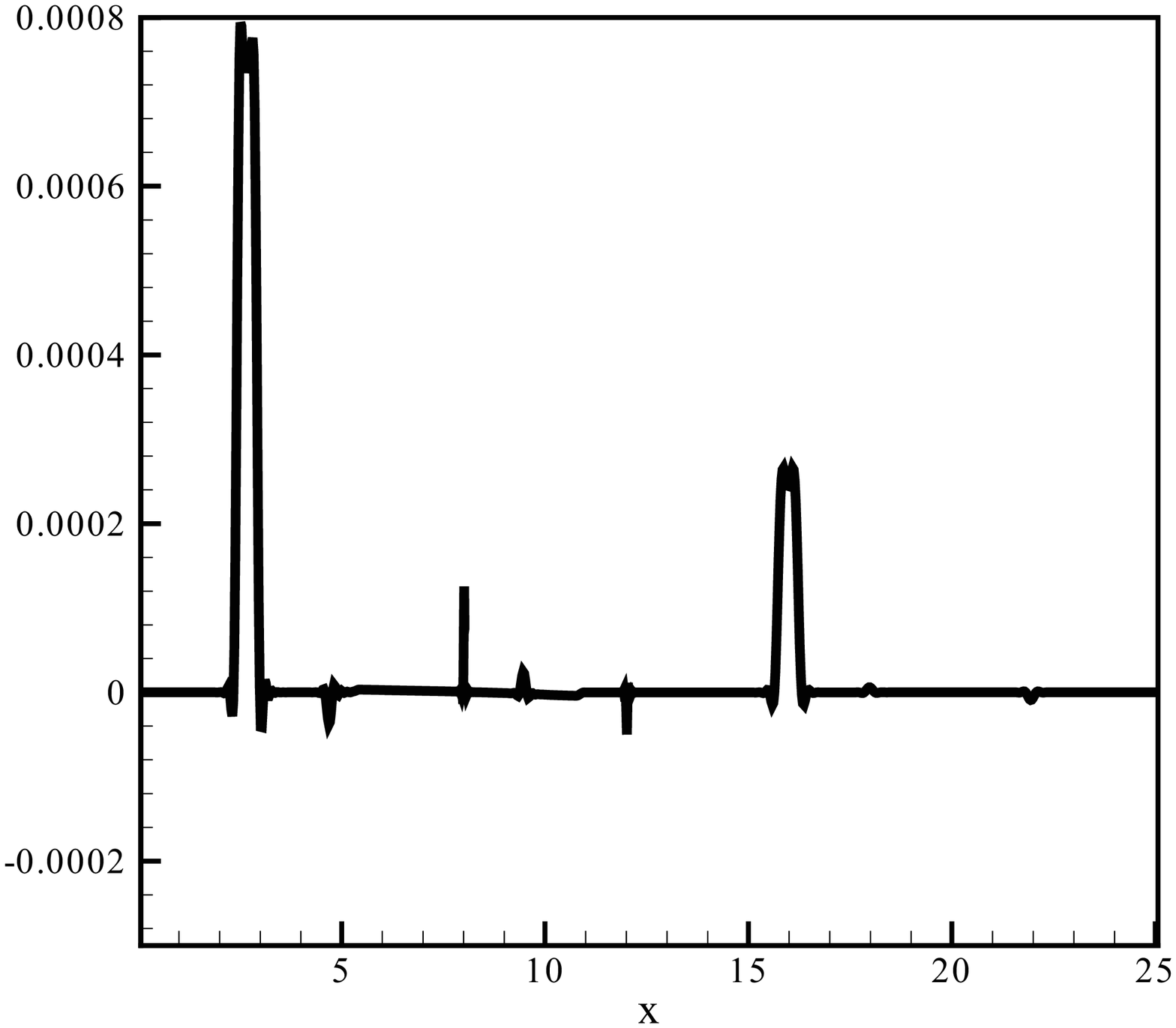}
\hspace{0.05in}
\includegraphics[width=3.2in]{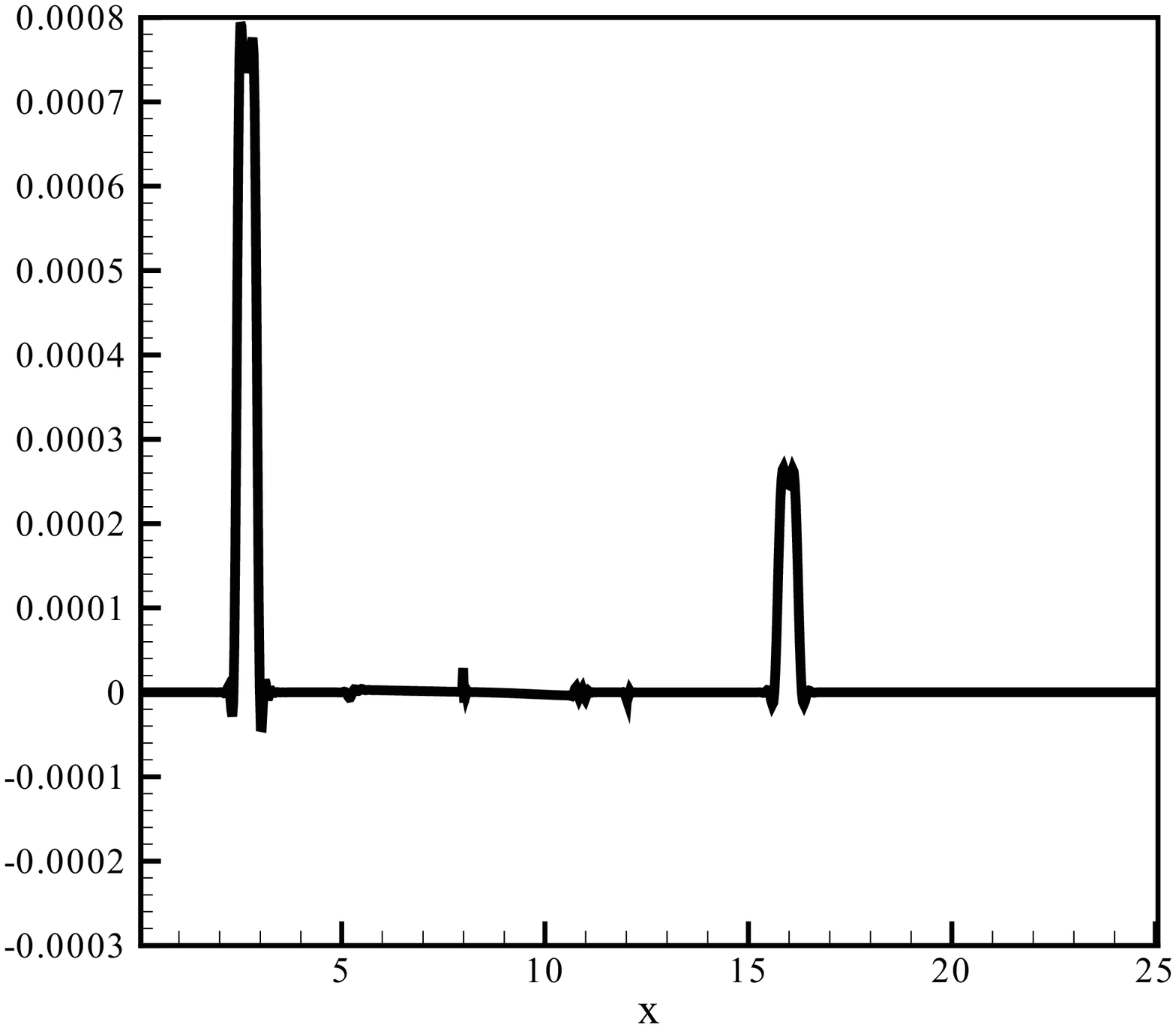}
} \vspace{-0.15in} \caption{Same as in Figure \ref{f11}, but 1000
uniform cells are employed.} \label{f13}
\end{figure}

\begin{figure}
\centerline{
\includegraphics[width=3.2in]{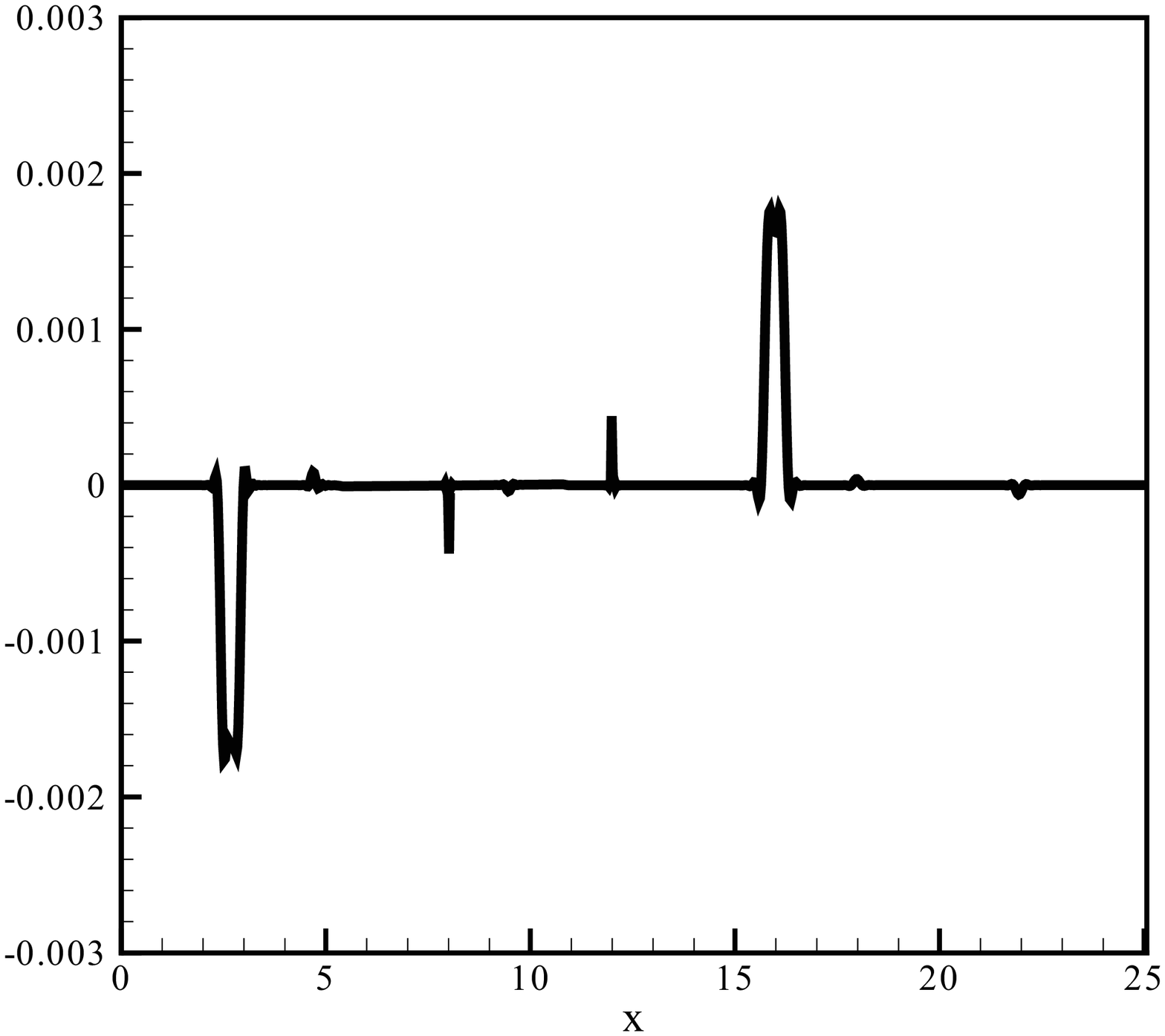}
\hspace{0.05in}
\includegraphics[width=3.2in]{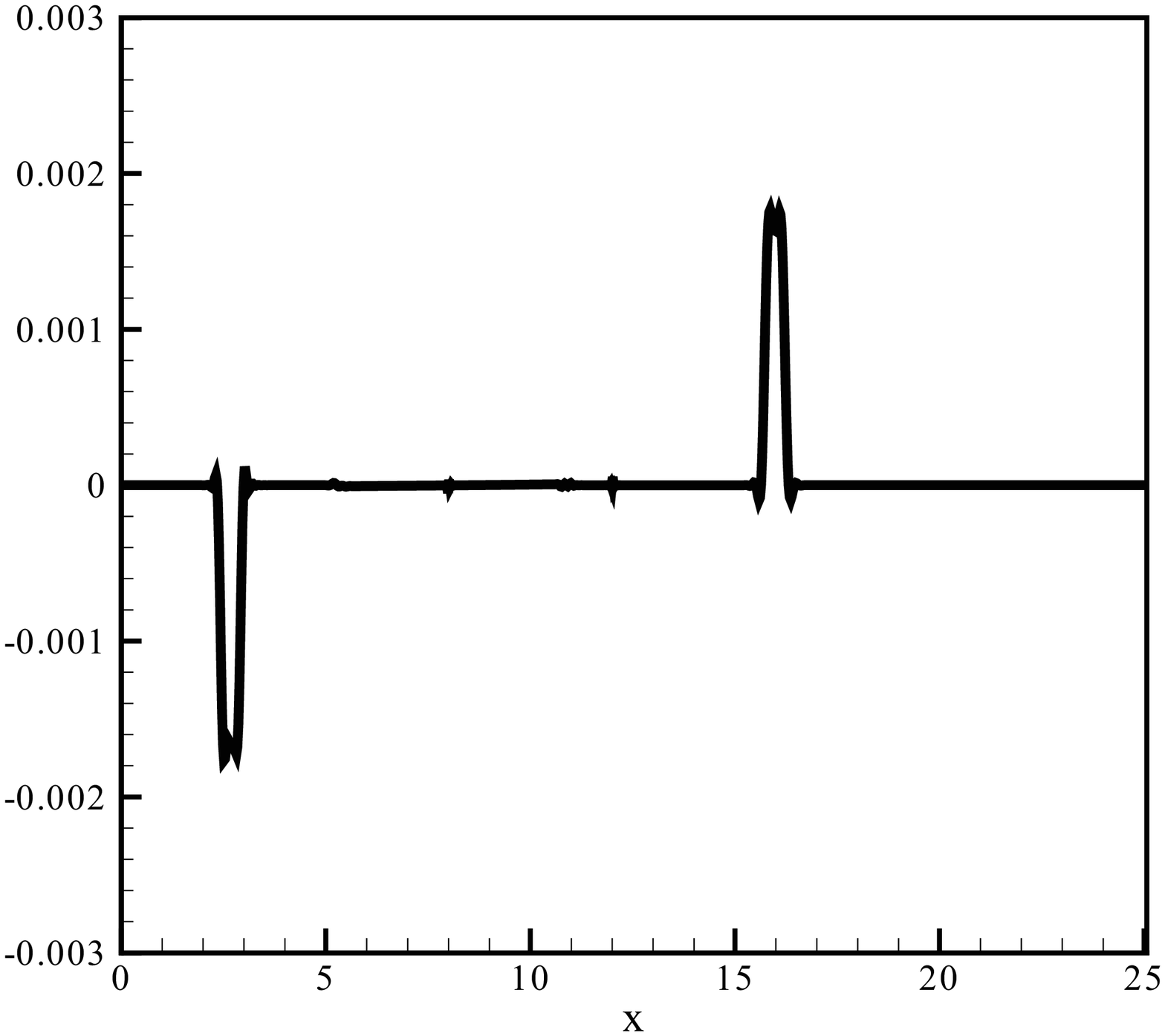}
} \vspace{-0.15in} \caption{Same as in Figure \ref{f21}, but 1000
uniform cells are employed.} \label{f23}
\end{figure}

We keep the main setup the same as in Section \ref{sec:n1}, and impose a
smaller perturbation of size $0.001$ to these steady states. To
save space, only the results of the subcritical flow
(\ref{eq:subcritical}) are shown. The differences between the
water height $h$ at that time and the background moving water
state, when 100 uniform cells are used, are plotted in Figure
\ref{f11}. We also show the differences of the momentum $hu$ in
Figure \ref{f21}. These figures clearly demonstrate that the still-water
well-balanced method is not capable of capturing such a small
perturbation on the coarse mesh, as large spurious oscillations are
observed. As we
refine the meshes to 200 cells, the results are shown in Figures
\ref{f12} and \ref{f22}, where spurious oscillations
for the still-water
well-balanced method have a
reduced magnitude, but its performance is still significantly
inferior to that of the moving-water well-balanced method.
The results with $1000$ cells are shown in Figures
\ref{f13} and \ref{f23}.  Even though the difference between the
two methods is now significantly reduced, the advantage of the
moving-water well-balanced method can still be observed on such
a refined mesh.

\medskip
{\bf Acknowledgement.} The first author is 
a contractor [UT-Battelle, manager of Oak Ridge National
Laboratory] of the U.S. Government under Contract No.
DE-AC05-00OR22725. Accordingly, the U.S. Government retains a
non-exclusive, royalty-free license to publish or reproduce the
published form of this contribution, or allow others to do so, for
U.S. Government purposes.

\end{document}